
\input psfig
\input amssym.def
\input amssym
\magnification=1100
\baselineskip = 0.25truein
\lineskiplimit = 0.01truein
\lineskip = 0.01truein
\vsize = 8.5truein
\voffset = 0.2truein
\parskip = 0.10truein
\parindent = 0.3truein
\settabs 12 \columns
\hsize = 5.4truein
\hoffset = 0.4truein

\setbox\strutbox=\hbox{%
\vrule height .708\baselineskip
depth .292\baselineskip
width 0pt}
\font\caps=cmcsc10

\def\sqr#1#2{{\vcenter{\vbox{\hrule height.#2pt
\hbox{\vrule width.#2pt height#1pt \kern#1pt
\vrule width.#2pt}
\hrule height.#2pt}}}}
\def\square{\mathchoice\sqr46\sqr46\sqr{3.1}6\sqr{2.3}4}
\def\leaderfill{\leaders\hbox to 1em{\hss.\hss}\hfill}
\font\bigtenrm=cmr10 scaled 1400
\tenrm

\centerline{\bf {\bigtenrm CLASSIFICATION OF ALTERNATING KNOTS}}
\centerline{\bf {\bigtenrm WITH TUNNEL NUMBER ONE}}
\tenrm
\vskip 14pt
\centerline{MARC LACKENBY}
\vskip 18pt

\centerline{\caps 1. Introduction}

An alternating diagram encodes a lot of information about a knot.
For example, if an alternating knot is composite, this
is evident from the diagram [10]. Also, its genus ([3], [12]) and
its crossing number ([7], [13], [17]) can be read off directly.
In this paper, we apply this principle to
alternating knots with tunnel number one.
Recall that a knot $K$ has {\sl tunnel number one}
if it has an {\sl unknotting tunnel}, which is defined to be
an arc $t$ properly embedded in the knot exterior
such that $S^3 - {\rm int}({\cal N}(K \cup t))$ is a handlebody.
It is in general a very difficult problem to determine
whether a given knot has tunnel number one, and
if it has, to determine all its unknotting tunnels.
In this paper, we give a complete classification of
alternating knots with tunnel number one, and all their
unknotting tunnels, up to an ambient isotopy of the
knot exterior.

\noindent {\bf Theorem 1.} {\sl Let $D$ be a reduced alternating diagram
for a knot $K$. Then $K$ has tunnel number one with an
unknotting tunnel $t$, if and only if $D$ (or its reflection)
and an unknotting tunnel isotopic to $t$ are 
as shown in Figure 1.
Thus the alternating knots with tunnel number one
are precisely the two-bridge knots and the Montesinos
knots $(e; p/q, \pm 1/2, p'/q')$ where $q$ and $q'$ are odd.}

\vskip 18pt
\centerline{\psfig{figure=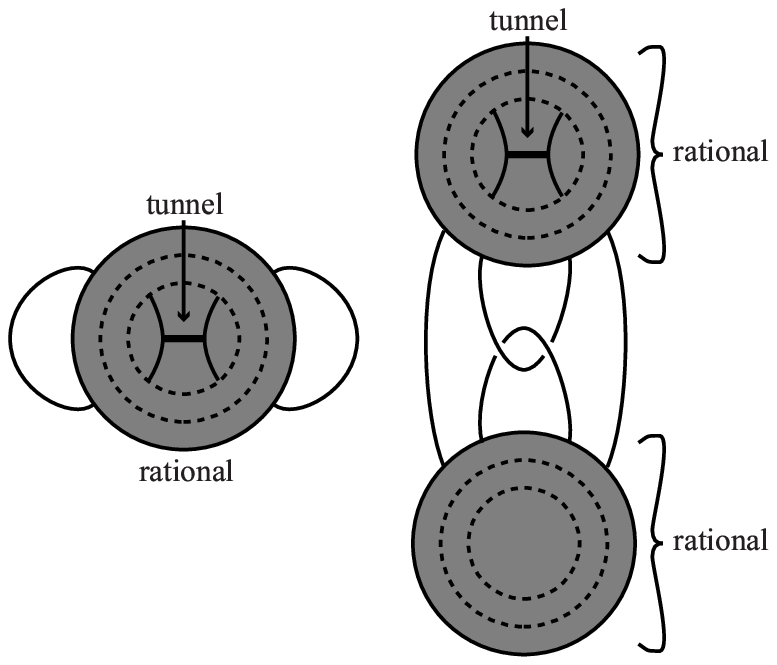,width=2.3in}}
\vskip 18pt
\centerline{Figure 1.}

For an explanation of the Montesinos knot terminology,
see [1]. The grey discs in Figure 1 denote alternating diagrams
of rational tangles with no nugatory crossings. 
It is proved in [18] (see the comments after Corollary 3.2 of [18])
that such a diagram is constructed
by starting with a diagram of a 2-string tangle containing
no crossings and then surrounding this diagram by
annular diagrams, each annulus containing four arcs
joining distinct boundary components, and each annular
diagram having a single crossing. (See Figure 2.)
The boundaries of these annuli
are denoted schematically by dashed circles within the grey disc.
Of course, the signs of these crossings are chosen
so that the resulting diagram is alternating. 

\vskip 18pt
\centerline{\psfig{figure=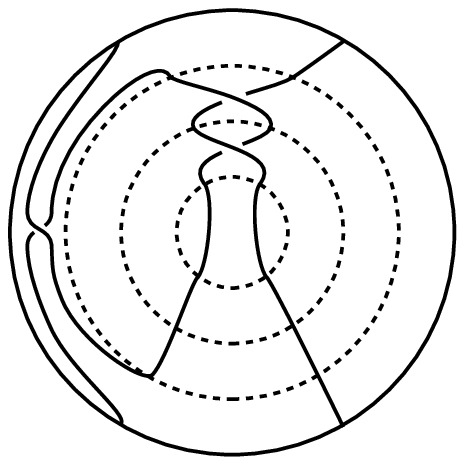,width=1.5in}}
\vskip 18pt
\centerline{Figure 2.}

Theorem 1 settles a conjecture of Sakuma, who
proposed that, when $D$ is a reduced alternating
diagram of a tunnel number one knot $K$, then some
unknotting tunnel is a vertical arc at some crossing of $D$.
For, the tunnels in Figure 1 can clearly be ambient
isotoped to be vertical at some crossing, unless the rational
tangle containing the tunnel has no crossings. But, in
this case, the diagram can be decomposed as in the left of Figure 1.
In particular, the knot is a 2-bridge knot.
Hence the knot has a vertical unknotting tunnel.

A given alternating diagram $D$ may sometimes be decomposed into the
tangle systems shown in Figure 1 in several distinct ways.
Hence, the knot may have several unknotting tunnels.
For example, from Theorem 1, it is not hard to deduce
Kobayashi's result [8], classifying all unknotting tunnels for a 2-bridge
knot into at most six isotopy classes.

We now explain why the knots in Figure 1 have tunnel number
one. In the left-hand diagram of Figure 1, contract the
tunnel to a point, resulting in a graph $G$.
It is clear that, by an ambient isotopy of $G$,
we may undo the crossings starting with the innermost annulus
and working out. Hence, the exterior of $G$ is a handlebody,
as required. 

\vskip 18pt
\centerline{\psfig{figure=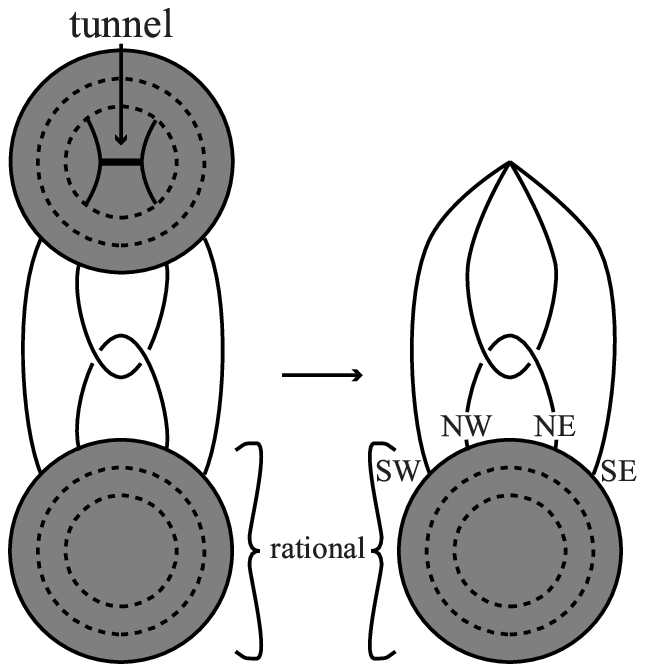,width=2in}}
\vskip 6pt
\centerline{Figure 3.}

We may perform a similar procedure in the right-hand diagram
of Figure 1, resulting in Figure 3. 
Pick the crossing in the outermost annulus
of the lower rational tangle. If it connects the strings
emanating from points SW and SE, we may remove it.
If it connects the strings emanating from NW and NE, we may
flype the rational tangle so that instead the crossing
lies between SW and SE, and then remove the crossing. 
If the crossing lies between points
NW and SW (or NE and SE) we may slide the graph as shown in
Figure 4 without altering the exterior, to change the
crossing. This procedure does not alter the way that
the two strings of the tangle join the four boundary
points. If the original diagram had more than one crossing,
then we consider the possibilities for the next annulus
in the inwards direction. By flyping if necessary, we may
assume that its crossing joins SW and SE, or SW and NW.
In both these cases, the tangle has an alternating diagram with
fewer crossings. Hence, inductively, we reduce to the case where
the original tangle has at most one crossing. 
Since $K$ is a knot, rather than a link, the strings of
the tangle run from NW to either SW or SE, and from NE to
either SE or SW. So, it is clear that the exterior of this graph is
a handlebody.

\vskip 18pt
\centerline{\psfig{figure=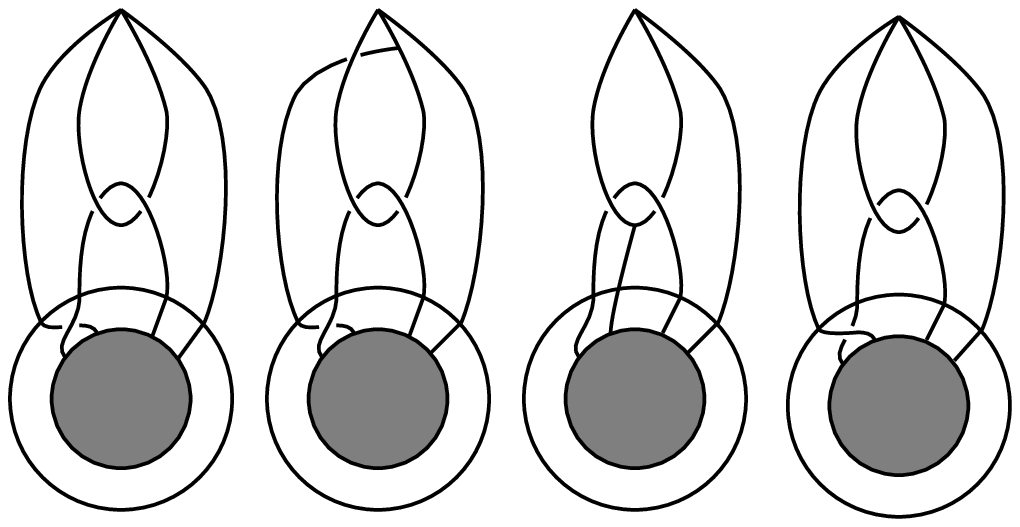,width=4in}}
\vskip 18pt
\centerline{Figure 4.}

Theorem 1 has been proved by Shimokawa [15] in the case where
the unknotting tunnel is isotopic to an
embedded arc in some region of the diagram. 
A large part of this paper is a proof that this must
always be the case.

\noindent {\bf Theorem 2.} {\sl Let $D$ be a reduced alternating diagram
for a knot $K$. Then any unknotting tunnel for $K$ is
isotopic to an unknotting tunnel that is an embedded arc 
in some region of the diagram.}

This leaves open the question of which planar arcs
are actually unknotting tunnels. This was answered by
Shimokawa [15], but we give a new proof of his result.
If one retracts a planar arc to a point, one does not change the exterior, 
but the resulting graph now has a single 4-valent vertex. We therefore
introduce the following definition.

\noindent {\bf Definition.} A diagram of an embedded
4-valent graph in $S^3$ is {\sl alternating} if there is
a way of modifying the vertices to crossings
so that the result is an alternating link diagram.
The diagram is {\sl reduced} if it is not
of the form shown in Figure 5.

\vskip 18pt
\centerline{\psfig{figure=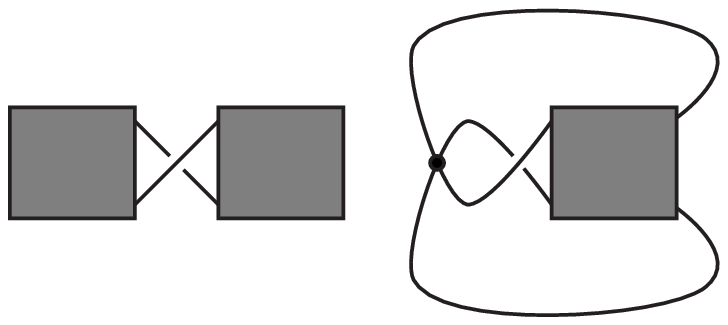}}
\vskip 18pt
\centerline{Figure 5.}

Morally, one perhaps should also consider
a diagram as in Figure 6, where each grey box
contains at least one crossing, as not reduced.
However, we will not adopt this convention in
this paper.

\vskip 18pt
\centerline{\psfig{figure=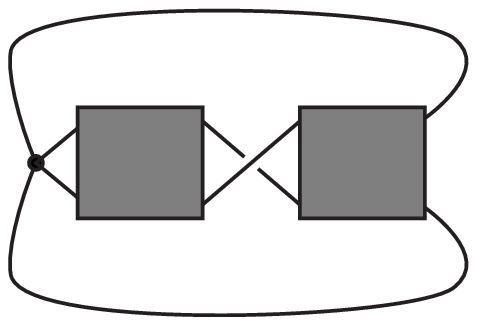}}
\vskip 18pt
\centerline{Figure 6.}

We prove the following result which can be used to
determine whether the planar arc is an unknotting tunnel.

\noindent {\bf Theorem 3.} {\sl Let $D$ be a reduced
alternating diagram of a graph $G$ with a single
vertex. Then the exterior of $G$ is a handlebody
if and only if $D$ is one of the diagrams
shown in Figure 7.}

\vskip 18pt
\centerline{\psfig{figure=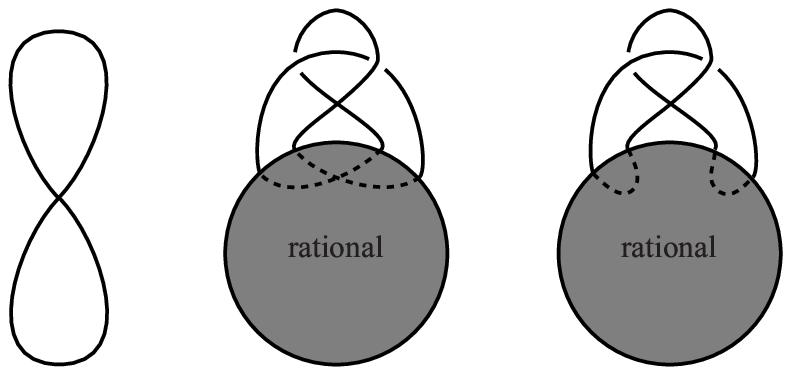}}
\vskip 18pt
\centerline{Figure 7.}

The dotted arcs in a grey disc denote the way that
the strings of the tangle join the four points
on the boundary. Note that the middle and right-hand alternating graphs 
in Figure 7 are ambient isotopic to that shown in Figure 3.

Theorem 1 is a straightforward corollary of Theorems
2 and 3. Given a reduced alternating diagram $D$ of a
knot $K$ with an unknotting tunnel $t$, we use Theorem 2 to 
establish that $t$ is isotopic to a planar arc.
Contract this tunnel to a point
to form a graph $G$ with an alternating diagram. 
Alter this diagram until it is reduced,
by removing crossings adjacent to the vertex.
Theorem 3 implies that it is of the form shown in
Figure 7. Hence, $D$ and the unknotting tunnel are as shown 
in Figure 1. This proves Theorem 1.

The proof of Theorems 2 and 3 is a largely amalgamation of ideas
by Rubinstein and Menasco. We now give
a brief outline of the main arguments.
It is a well-known fact that an unknotting tunnel
for a knot is determined by its associated Heegaard
surface, up to an ambient isotopy of the
knot exterior. This follows from
the fact that, in a compression body $C$ with
$\partial_- C$ a torus and $\partial_+ C$ a genus
two surface, $\partial C$ has a unique non-separating
compression disc up to ambient isotopy. Hence, we will focus
on the Heegaard surface.
Rubinstein showed that, given any triangulation
of a compact orientable irreducible 3-manifold $M$, a strongly irreducible
Heegaard surface for $M$ can be ambient isotoped into almost
normal form [14]. An alternating knot complement inherits an ideal polyhedral
structure from its diagram [9]. We therefore in \S2
develop a notion of almost normal surfaces in such an ideal
polyhedral decomposition of a 3-manifold. A genus two 
Heegaard surface $F$ for a non-trivial knot must be strongly irreducible, and
so can be ambient isotoped into almost normal form. 
Once in this form, $F$ intersects the plane
of the diagram in a way rather similar to
the surfaces studied by Menasco [10]. In \S3,
we recall Menasco's techniques. 

It would seem logical to prove Theorem 2 before Theorem 3,
but in fact the latter result is necessary in the
proof of the former. Therefore in \S4, we prove Theorem 3.
The hypothesis that the graph's exterior is a handlebody
is used to establish the existence of a compressing
disc in normal form. This is then used to show that
$G$ must be as described in the theorem.

In \S5, we adapt Menasco's
arguments to restrict the possibilities for the
Heegaard surface $F$. We show that $F$ is obtained
from a standardly embedded 4-times punctured sphere by attaching tubes
that run along the knot. The sphere divides the diagram
into two tangles. We analyse the
possibilities for these tangles, using Theorem 3,
and prove that the unknotting tunnel 
is isotopic to a planar arc in one of them.
This will prove Theorem 2 and hence Theorem 1.

\vfill\eject
\centerline{\caps 2. Heegaard surfaces in ideal polyhedral decompositions}
\vskip 6pt

A reduced alternating diagram of knot determines an ideal polyhedral
decomposition of the knot complement. Hence, in this section,
we develop an almost normal surface theory for ideal
polyhedral decompositions of 3-manifolds and establish that, under
certain conditions,
a strongly irreducible Heegaard surface can be ambient isotoped
into almost normal form.

\noindent {\bf Definition.} A {\sl polyhedron} is a 3-ball
with a non-empty connected graph in its boundary, the graph having no
edge loops.
An {\sl ideal polyhedron} is a polyhedron with its vertices removed.
An {\sl ideal polyhedral decomposition} of a 
3-manifold $M$ is a way of constructing $M - \partial M$ as a 
union of ideal polyhedra with their faces identified in pairs.

Note that polyhedra cannot necessarily be realised geometrically
in Euclidean space with straight edges and faces. In particular,
it is possible for a face in a polyhedron to have only
two edges in its boundary. Such faces are known as {\sl bigons}.

There is a reasonably well-established theory of normal
surfaces in ideal polyhedra. Associated with an
ideal polyhedral decomposition of a 3-manifold $M$,
there is a dual handle decomposition of $M$, where $i$-handles
($0 \leq i \leq 2$) arise from $(3 - i)$-cells in the
ideal polyhedra. In [6], Jaco and Oertel gave a 
definition of a closed normal surface in such
a handle decomposition. Here we give the dual version.

\noindent {\bf Definition.} A disc properly embedded in
an ideal polyhedron $P$ is {\sl normal} if
\item{$\bullet$} it is in general position with respect
to the boundary graph of $P$, 
\item{$\bullet$} it intersects each edge of $P$ at most once, and
\item{$\bullet$} its boundary does not lie wholly in some face of $P$.

\noindent A closed properly embedded surface in a 3-manifold
$M$ is in {\sl normal form} with respect to an ideal polyhedral 
decomposition of $M$ if it intersects each ideal polyhedron
in a (possibly empty) collection of normal discs.

Note that a normal surface intersects each face in {\sl normal}
arcs, which means that each arc is not parallel to a sub-arc of an edge.

There is also notion of normality for properly embedded
surfaces with non-empty boundary. We will come to this
at the end of \S3.

The following is an elementary fact about normal surfaces.
The proof is an easy generalisation of the case where each
polyhedron is a tetrahedron, which can be found in [19].

\noindent {\bf Lemma 4.} {\sl Fix an ideal polyhedral decomposition
of a 3-manifold in which each face is a triangle or a bigon.
Then a closed normal surface is incompressible in the complement
of the 1-skeleton.}

Therefore, given any ideal polyhedral decomposition of a 3-manifold,
we will always subdivide its faces into bigons and triangles,
by possibly introducing new edges, but adding no vertices.

It is a well-known result [6] that
any closed properly embedded incompressible surface with no 2-sphere
components in a compact orientable irreducible 3-manifold $M$
can be ambient isotoped into normal form with respect to
some fixed triangulation. Stocking [16], building on ideas of Rubinstein [14]
and Thompson [19], proved that any strongly irreducible
Heegaard surface in $M$  can be ambient isotoped into `almost normal' form.
Recall that a closed surface is {\sl almost normal}
with respect to a triangulation of $M$ if its
intersection with each tetrahedron is a collection
of squares and triangles, except in precisely one
tetrahedron, where precisely one component of 
intersection between the surface and the tetrahedron
is either an `octagon' or an annulus obtained
from two triangles or squares by attaching
an unknotted tube. These two possibilities are shown
in Figure 8. 

\vskip 18pt
\centerline{\psfig{figure=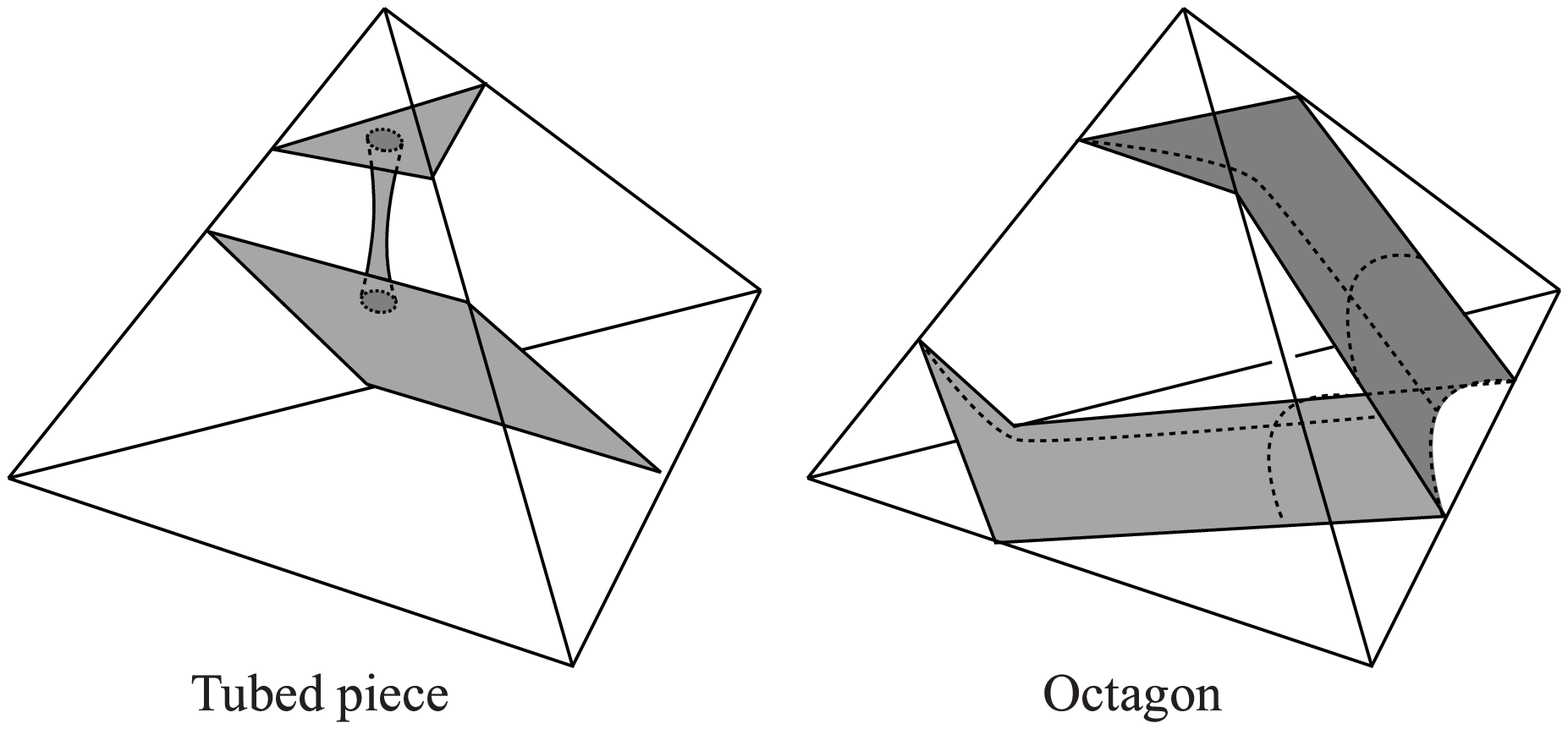,width=3.5in}}
\vskip 18pt
\centerline{Figure 8.}

Stocking's result holds true even when $M$ has
non-empty boundary: in this case, we consider the usual definition [2]
of a Heegaard splitting for $M$ as a decomposition into compression bodies.
Recall that a {\sl compression body} $C$ is a connected compact orientable
3-manifold that is either a handlebody or is
obtained from $S \times [0,1]$, where $S$ a 
closed (possibly disconnected) surface, by attaching
1-handles to $S \times \lbrace 1 \rbrace$.
Then we let $\partial_-C$ be $S \times \lbrace 0 \rbrace$,
or the empty set when $C$ is a handlebody, and we
let $\partial_+ C = \partial C - \partial_- C$.
A {\sl Heegaard splitting} for a compact orientable 3-manifold $M$ is a
decomposition of $M$ into two compression bodies
glued along their positive boundaries.

One feature of octagonal and tubed pieces is 
that the resulting almost normal surface $F$ has an
{\sl edge compression disc}, namely a disc
$D$ embedded in a tetrahedron $\Delta$
such that $D \cap \partial \Delta$ lies in $\partial D$ and is a sub-arc
of an edge, and $D \cap F$ is the remainder of $\partial D$.
There is one aspect of polyhedral decompositions that makes them
a little more complicated than triangulations. If a disc
properly embedded in a tetrahedron intersects each face
in a collection of normal arcs
and has an edge compression disc on one side, then it
has an edge compression disc on the other side also. This
need not be true in more general polyhedra. Therefore
we introduce the following definition.

\noindent {\bf Definition.} Let $F$ be a closed two-sided surface
properly embedded in a compact orientable 3-manifold
with an ideal polyhedral decomposition. Then $F$ is {\sl normal
to one side} if
\item{(i)} its intersection with any face is a collection
of normal arcs,
\item{(ii)} its intersection with any ideal polyhedron
is a collection of discs, and
\item{(iii)} all edge compression discs for any of these
discs emanate from the same side of $F$.

Note that a surface that is normal to one side may in fact
be a normal surface. However, we will often refer to `the' normal
side of the surface, with the understanding that if the surface is
normal, then this means either side. We now introduce a generalisation of almost
normality to surfaces in ideal polyhedral decompositions.

\vfill\eject
\noindent {\bf Definition.} Let $P$ be an ideal polyhedron in which each face 
is a triangle or bigon. An {\sl almost normal disc} in $P$ is
a properly embedded disc $D$ having the following properties:
\item{(i)} its intersection with each face is
a collection of normal arcs, 
\item{(ii)} in each component of ${\rm cl}(P - D)$
there is an edge compression disc, and
\item{(iii)} any two edge compression discs, one on each side of $D$,
must intersect in the interior of $P$.

\noindent Let $M$ be a compact orientable 3-manifold with an
ideal polyhedral decomposition in which each face is a
triangle or bigon. A closed properly embedded surface in $M$
is {\sl almost normal} if either
\item{(a)} its intersection with
each ideal polyhedron is a collection of normal discs,
except in precisely one ideal polyhedron where it
has precisely one almost normal disc, together possibly 
with some normal discs, or
\item{(b)} it is obtained from a normal surface by attaching 
a tube that lies in a single ideal polyhedron, that runs
parallel to an edge and that joins distinct normal discs.

An example of an almost normal disc is shown in Figure 9.

\vskip 18pt
\centerline{\psfig{figure=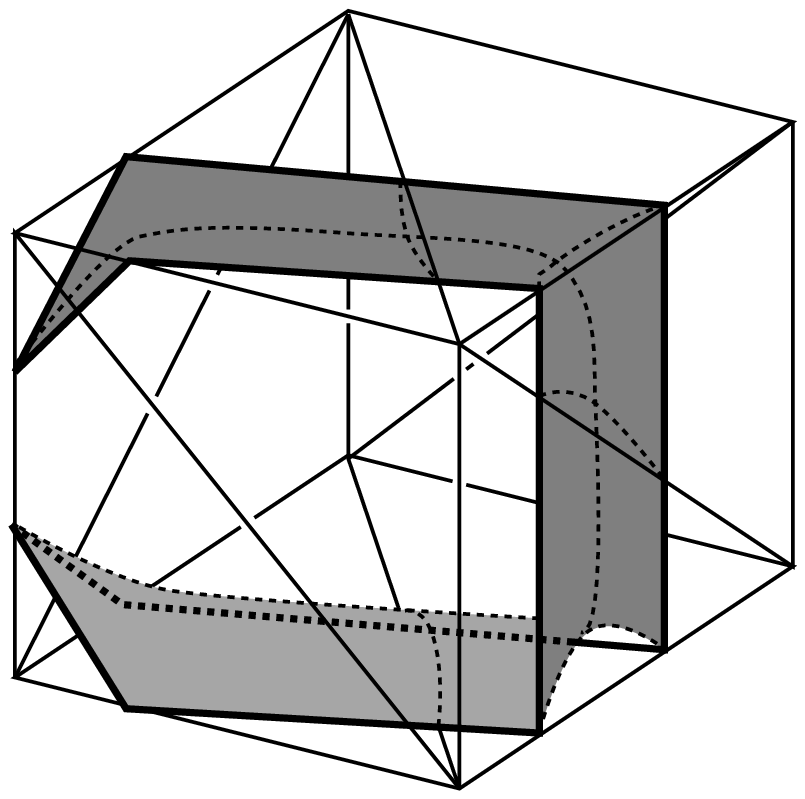,width=2in}}
\vskip 18pt
\centerline{Figure 9.}
\vfill\eject

\noindent {\bf Theorem 5.} {\sl Let $F$ be a strongly irreducible
Heegaard surface for a compact orientable irreducible 3-manifold $M$.
Fix an ideal polyhedral decomposition of $M$ in which
each face is a bigon or triangle. Suppose that $M$ contains
no 2-spheres that are normal to one side. Then
there is an ambient isotopy taking $F$ into almost
normal form.}

The proof of this theorem follows Stocking's argument [16]
very closely. We therefore will only sketch the main
outline and will refer the reader to [16] for further details.

The hypothesis that $M$ contains no 2-spheres that are
normal to one side is an unfortunate one. Although this
theorem is sufficient for our purposes, it is possible
that this assumption can be dropped, with some substantial
work. Normal 2-spheres caused significant complications
in Stocking's proof, and the possibility of 2-spheres
that are normal to only one side causes even more difficulties.

We will need the following lemma, which is a translation
of Lemma 1 in [16] to the polyhedral setting. However,
the proof in [16] involves a case-by-case analysis that
is not suitable here. We therefore give an alternative
proof.

\noindent {\bf Lemma 6.} {\sl Let $S$ be a closed separating properly embedded
surface in a compact orientable irreducible 3-manifold $M$
with an ideal polyhedral decomposition in which each face
is a triangle or bigon. Suppose that $S$ is
almost normal or normal to one side. Suppose also that
$S$ is incompressible into one component $I$ of $M - 
{\rm int}({\cal N}(S))$, and that $S$ has an edge compression
in $I$. Then $S$ is ambient isotopic to a surface
in $I$, each component of which is either normal or a
2-sphere lying entirely in a single ideal polyhedron.}

\noindent {\sl Proof.} 
We will perform a sequence of isotopies that move $S$ into $I$. At each stage,
we will denote the new surface by $S$ and the new
component of $M - {\rm int}({\cal N}(S))$ into which 
$S$ is incompressible by $I$. We are assuming that $S$ has an
edge compression disc in $I$. This specifies 
an ambient isotopy that reduces the weight of surface. Perform
this ambient isotopy into $I$. The result is a surface 
that need not be normal, but which has the following properties:
\item{(i)} for any component $B$ of intersection between $I$ and
a polyhedron, the intersection between $B$ and any
face  has at most one component;
\item{(ii)} if an arc of intersection between
$S$ and some face has endpoints on the same edge,
and $D$ is the subdisc of the face that it separates off,
then the interior of $D$ is disjoint from $S$
and is a subset of $I$;
\item{(iii)} if $C$ is a simple closed curve of intersection
between $S$ and some face, and $D$ is the subdisc 
of the face that it separates off,
then the interior of $D$ is disjoint from $S$
and is a subset of $I$;
\item{(iv)} any component of intersection between $S$ and a
polyhedron $P$ is incompressible in the $M-I$ direction;
\item{(v)} for any component $B$ of intersection between $I$
and a polyhedron $P$, $\partial B \cap P$ is boundary-parallel in $P$.

Note that (i) need not apply to $S$ before the isotopy. But
it does afterwards. Otherwise, we may find an edge compression
disc for the original $S$ in $M - I$ with the property
that it and the edge compression disc in $I$ do not intersect
away from the 2-skeleton of $M$. This contradicts the assumption that $S$ is
almost normal or normal to one side.

Suppose that in some polyhedron $P$, $S \cap P$ is not a
collection of discs and 2-spheres. Then $S \cap P$ has a compression disc $D$
in $P$. This disc lies in $I$, by (iv). Hence, $\partial D$
bounds a disc $D'$ in $S$. Ambient isotope $D'$ onto $D$.

Note that if $S$ intersects some face in a simple closed curve,
but $S \cap P$ does not admit a compression in any polyhedron $P$, 
then this is the only intersection between this component of $S$
and the 2-skeleton. By (iii), this 2-sphere bounds a 3-ball in $I$,
and we may ambient isotope it into a single ideal polyhedron.

Suppose that, for some polyhedron $P$ and some
component $D$ of $S \cap P$, $D$ intersects some
edge $e$ more than once. We may assume that there
are two adjacent points of intersection between
$D$ and $e$. Hence, there is an edge compression disc.
This edge compression disc must lie in $I$, otherwise
we contradict (i) or (ii). Hence, we use this to perform
an isotopy into $I$, reducing the weight of the surface.

After each of these isotopies, properties (i) - (v) still hold.
Eventually, this process must terminate in the
required surface. $\square$

\noindent {\bf Corollary 7.} {\sl If $M$ contains no normal
2-spheres, then it contains no almost normal 2-spheres.}

\noindent {\sl Proof.} Let $S$ be an almost normal 2-sphere,
and let $M_1$ and $M_2$ be the two components of $M - {\rm int}
({\cal N}(S))$. By Lemma 6, we may ambient isotope
$S$ into each $M_i$ until it is normal or disjoint from
the 2-skeleton. The former case is impossible, by hypothesis. Thus,
$S$ must bound a 3-ball in both sides. Therefore, $M$ is
the 3-sphere, which is closed, but closed 3-manifolds do not
have an ideal polyhedral decomposition. $\square$

We will prove Theorem 5 by induction. At each stage, we will
consider a connected compact 3-manifold $M_i$ embedded in $M$
such that
\item{(i)} $\partial M_i$ is normal in $M$, and
\item{(ii)} $F$ lies in $M_i$ and is a Heegaard surface for $M_i$.

\noindent Initially, $M_1 = M$. Note that $M_i$ inherits an ideal
polyhedral decomposition from that of $M$, by taking the
intersection between $M_i$ and the ideal polyhedra of $M$,
and then removing $\partial M_i$. The hypothesis that
$\partial M_i$ is normal in $M$ guarantees that surfaces
that are normal, normal to one side or almost normal in
$M_i$ have the same property in $M$.

At each stage of the induction, we
either deliver the Heegaard surface $F$ in
almost normal form, or we construct a new embedded
3-manifold $M_{i+1}$ contained in $M_i$, that
satisfies (i) and (ii) above. The boundary of $M_{i+1}$
will not be parallel in the ideal polyhedral decomposition
to that of $M_i$. A straightforward modification of the standard
argument due to Kneser gives that
there is an upper bound on the number of
such surfaces in $M$ [5]. So, we eventually obtain
$F$ in almost normal form. We now give the main steps of the argument. 

The Heegaard splitting for $M_i$ determines a singular
foliation, as follows. In the case where a compression body $C$ is
a handlebody, the singular set in $C$
is a graph onto which $C$ collapses; otherwise it is
the cores of the 1-handles.
The complement of the singular set and $\partial_- C$ is given a product
foliation $F \times (0,1)$. A small isotopy guarantees
that the 1-skeleton $\Delta_1$ of the ideal polyhedral
decomposition is disjoint from the singular set. We then
place $\Delta_1$ in thin position with respect to $F \times (0,1)$.
The following proposition is the key step.

\noindent {\bf Proposition 8.} {\sl There is a non-singular leaf $F$ of
the foliation, which either is almost normal or 
can be compressed on one side in the complement of
the 1-skeleton to a (possibly disconnected) almost normal surface 
$\overline F$. In the latter case, the incompressible
side of $\overline F$ has an edge compression disc.}

\noindent {\sl Proof.}
There are two cases to consider: when there is a thick region
and when there is not. Consider first the case where there is no thick region.
Then we take $F$ to be a leaf in the foliation having
the maximal number number of points of intersection with $\Delta_1$.
This is obtained from $\partial_- C$ for one of the
compression bodies $C$, by attaching tubes, which are the boundaries
of small regular neighbourhoods of the singular set.
The argument of Lemma 4 in [16] gives that, after possibly
compressing some of these tubes to one side, we obtain an almost normal surface
that contains a tubed component. If any compressions were used,
the resulting surface is incompressible to one side, as $F$ is
strongly irreducible [2]. The edge compression disc for the
tube lies on this side.
(Essentially, Lemma 4 in [16] barely uses that the
hypothesis that $M$ is triangulated; instead it is
an analysis of how the tubes lie in the ideal
polyhedral decomposition.) 

Consider now the case where there is a thick region.
Applying the proof of Claim 4.4 in [19], we can find a leaf $F$ of
the foliation in the thick region, which intersects each
face in a collection of normal arcs and simple closed curves.
It has an upper and a lower disc. The assumption that
this is thin position guarantees that any two such
discs must intersect at other than their endpoints.
We compress $F$, if necessary, to a (possibly disconnected) 
surface $\overline F$ which intersects each face
only in normal arcs, and which intersects each ideal
polyhedron in discs. 

We claim that $\overline F$ is normal to one side or
almost normal. If $\overline F$ has edge compression discs on 
at most one side, it is normal to one side.
If it has edge compression discs on both sides, they 
must lie in the same polyhedron, otherwise we contradict thin position. 
By the argument of Claims 4.1 to 4.3 in [19], at most one disc of $\overline F$
in this polyhedron can be non-normal, and
it must be almost normal. We indicate briefly
how this argument runs. Any edge compression
disc for any disc of $\overline F$ can be isotoped so that
its interior is disjoint from $\overline F$. For,
otherwise, $\overline F$ (and hence $F$) 
has a pair of nested upper and lower discs,
contradicting thin position. Thus, if 
$\overline F$ contains two non-normal discs,
we may find disjoint edge compression discs,
one emanating from each side of $\overline F$.
These form disjoint upper and lower discs for
$F$, which  is a contradiction. Similarly,
any two edge compression discs for a disc of
$\overline F$, one on each side of $\overline F$,
must intersect away from their boundaries.
This proves the claim.

Now, $F$ is obtained from $\overline F$ by attaching tubes.
We claim that they are not nested, and that their
meridian discs all lie on the same side of $F$.
For if they are nested, then we may pick a tube $T$,
with at least one tube running through it, but such
that all tubes $T_1, \dots T_n$ running through $T$ are innermost.
Let $D_i$ be a meridian disc for $T_i$, and let $D$
be a meridian disc for $T$. If $D_i$
is not a compression disc for $F$, then $\partial D_i$
bounds a disc in $F$. If this disc contains any tubes,
consider the simple closed curves forming the boundaries
of their meridian discs. Pass to an innermost such curve.
This bounds a disc in $F$ which forms part of a 2-sphere
component of $\overline F$. But we have made the assumption
that $M$ contains no 2-spheres that are normal to one side.
By Corollary 7, $M$ also contains no almost normal
2-spheres. Thus, no component of $\overline F$ is
a 2-sphere. Thus, each $D_i$ is a compression disc
for $F$. So, when we compress $F$ along these discs,
the resulting surface $F'$ is incompressible on
the side containing $D$, as $F$ is strongly irreducible.
Therefore, $\partial D$ bounds
a disc in $F'$. Again, this implies that $\overline F$
contains a 2-sphere component, which is a contradiction. 
Therefore, the tubes of $F$ are not nested. Their
meridian discs are all essential. Since $F$ is strongly
irreducible, they all lie on the same side of $F$.
This proves the claim.

We claim that $\overline F$ is almost normal, giving
the required surface. Suppose that, on the contrary, $\overline F$
is normal to one side. Now, $F$ has both upper and lower discs. 
By the argument in Claim 11 of [16], we may assume that
they are disjoint from the interiors of the tubes. Thus,
on the side of $\overline F$ to which the tubes
are not attached, there must be an edge compression disc,
making that side non-normal.
So, the tubes are attached to the normal side of $\overline F$.
Then, since the tubes of $F$ are not nested and all emanate
from the same normal side, we may apply
by the argument of Lemma 4 in [16] to deduce that
there is an edge compression disc for $F$ that runs
over one tube exactly once. This and the edge
compression disc on the non-normal side of $\overline F$
form upper and lower discs for $F$ that intersect
away from $\Delta_1$, which is a contradiction. $\square$

If $F$ is almost normal, Theorem 5 is proved.
If not, then by Proposition 8, it compresses on one side
to an almost normal surface $\overline F$.
Let $I$ be the incompressible side of $\overline F$.
The tubes of $F$ are all attached to this side, and
hence lie in $I$. Apply Lemma 6 to ambient isotope
$\overline F$ into $I$ to a normal surface. 
Let $M_{i+1}$ be the copy of $I$ after the ambient isotopy.
This has the required properties. Hence the proof of
Theorem 5 is complete.

\vskip 18pt
\centerline{\caps 3. Polyhedral decompositions of
alternating knot complements}
\vskip 6pt

It is well known that a reduced diagram of a knot
induces a decomposition of the knot complement into two
ideal polyhedra. The construction is due to Menasco [9].
We recall the main details now.
Suppose we are given a reduced knot diagram $D$ lying in a
2-sphere.  We embed this 2-sphere into $S^3$. The knot
lies in this 2-sphere, except near each crossing, where
it skirts above and below the diagram as two semi-circular
arcs (see Figure 10). These arcs lie on the boundary of a
2-sphere `bubble' that encloses each crossing. The
2-sphere containing the diagram decomposes each bubble into
two discs, an upper and lower hemisphere. The upper (respectively, lower)
hemispheres together with the remainder of the diagram 2-sphere is
a 2-sphere denoted $S^2_+$ (respectively, $S^2_-$).
We will use $S^2_\pm$ to denote $S^2_-$ or $S^2_+$.

\vskip 18pt
\centerline{\psfig{figure=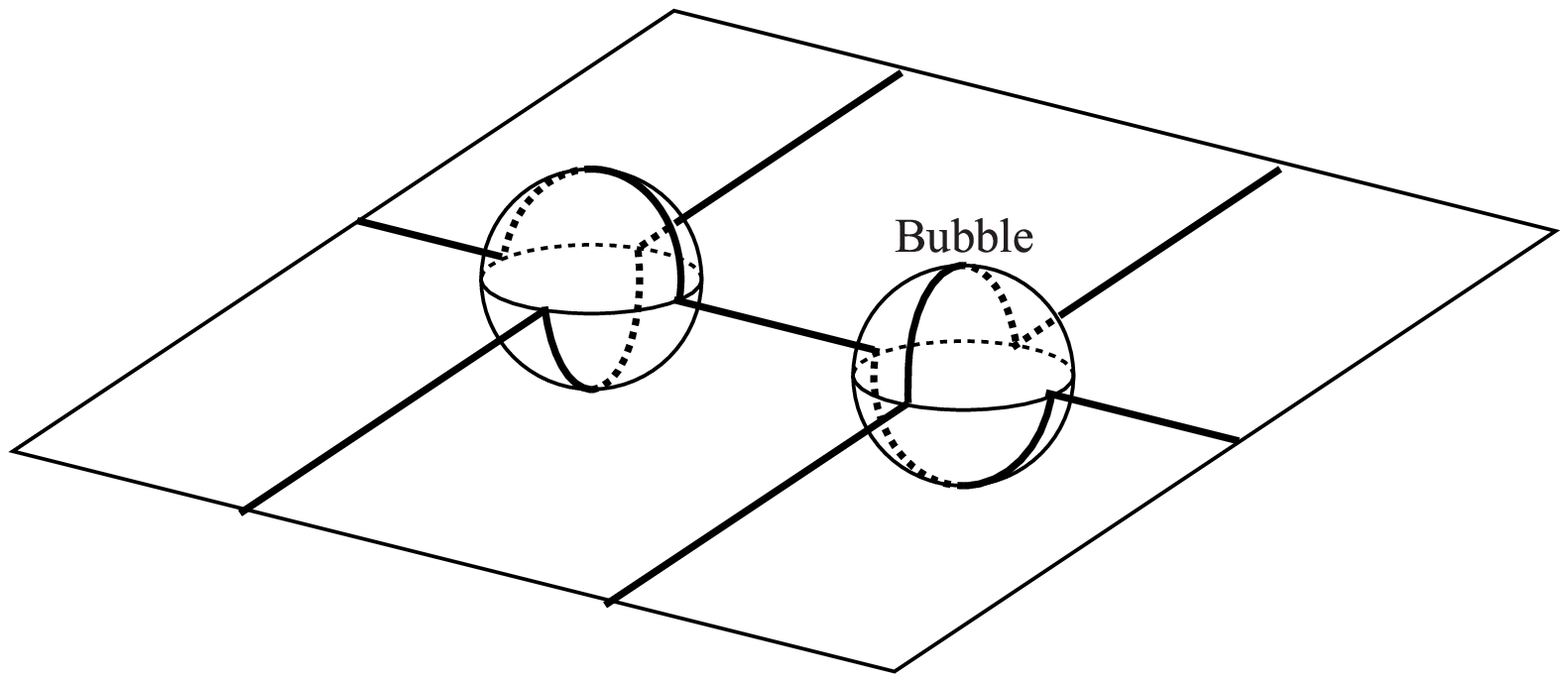,width=3in}}
\vskip 18pt
\centerline{Figure 10.}

At each crossing there is a vertical arc properly embedded in the knot complement.
These arcs form the edges of the ideal polyhedral decomposition.
There is one face for each region of the diagram. (See Figure 11.)
The complement of the edges and faces is two open 3-balls,
one above the diagram, one below. These open balls are
the interior of the two ideal polyhedra.

A closed normal or almost normal surface $F$ in this ideal polyhedral decomposition
can be visualised in quite a straightforward way. At each
point of intersection between $F$ and an edge, we insert
at a saddle in the relevant bubble (see Figure 12). The intersection between
$F$ and the faces of the polyhedra is a collection of arcs which
lie in the regions of the diagram. The condition that $F$ is normal
(or almost normal) and that the diagram is reduced
guarantees that no arc has endpoints in the same crossing.

\vskip 18pt
\centerline{\psfig{figure=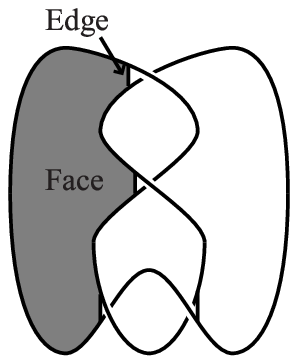}}
\vskip 18pt
\centerline{Figure 11.}

\vskip 18pt
\centerline{\psfig{figure=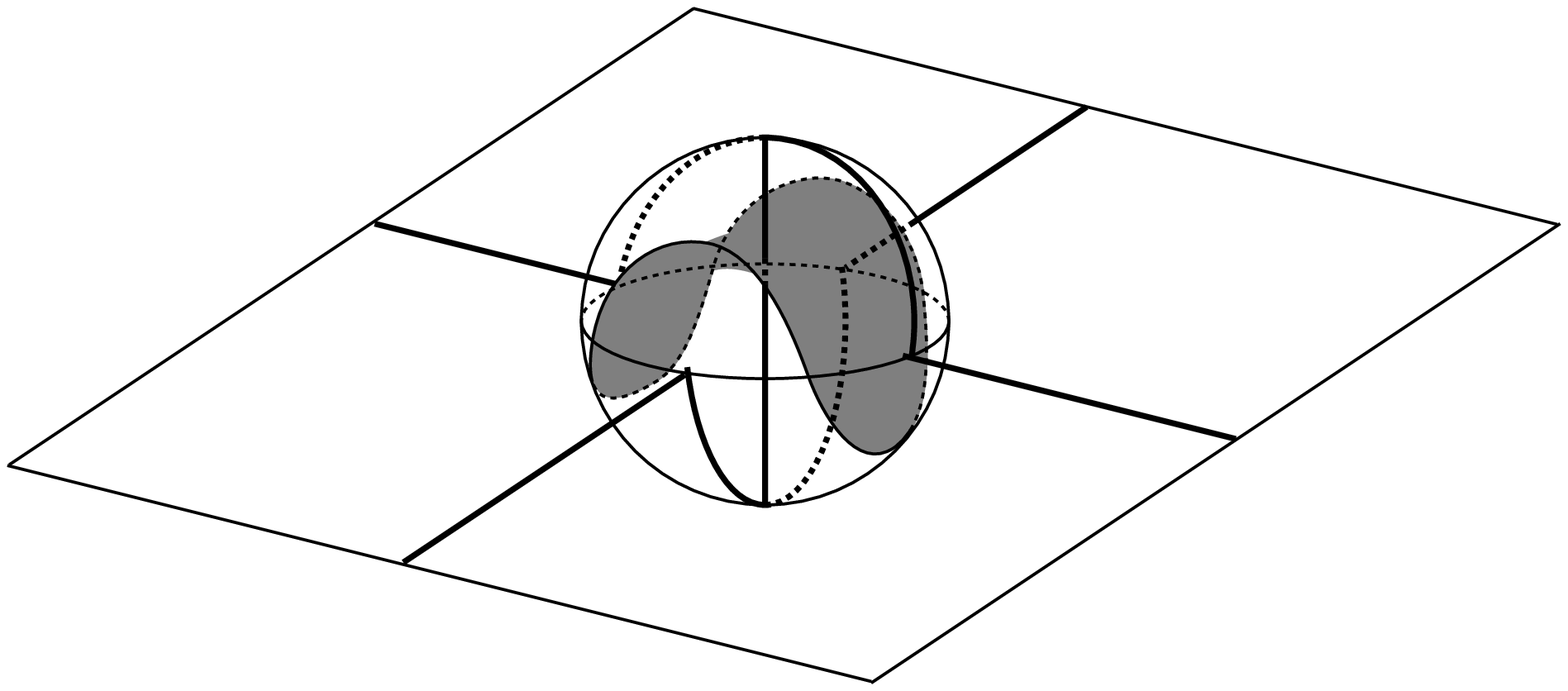,width=3in}}
\vskip 18pt
\centerline{Figure 12.}

We will need to consider surfaces that are more general
than normal surfaces. Let $F$ be a surface properly embedded
in the knot complement, such that any boundary component
of $F$ is meridional. We say that $F$ is {\sl standard} if it
intersects each bubble in a collection of saddles,
intersects each face in a collection of arcs
and intersects the two truncated polyhedra in
a collection of discs. Also, the boundary $C$ of each such disc
must satisfy the following conditions:
\item{(i)} no arc of intersection between $C$ 
and a face has endpoints lying in the same crossing;
\item{(ii)} if $C$ intersects the knot, it does so
transversely away from the bubbles;
\item{(iii)} $C$ intersects the knot projection in more than two
points.

A central part of Menasco's techniques [10] is to analyse how a
normal surface $F$ lies in an alternating knot complement. Under
many circumstances, he established the existence of a {\sl meridional
compression disc} which is defined to be an embedded disc
$R$ such that $R \cap F = \partial R$ and
$R \cap K$ is a single point in the interior of $R$.
The meridional compression disc that Menasco constructs is
{\sl diagrammatic} which means that its intersection with the
bubbles is a single disc in a single bubble, and the remainder
of the disc is disjoint from the plane of the diagram.

We say that a 2-sphere is {\sl trivial} if it is disjoint from the bubbles
and it intersects the plane of the diagram in a single
simple closed curve.

The following is a generalisation of Menasco's results
to standard surfaces.

\noindent {\bf Lemma 9.} {\sl Let $D$ be a prime alternating
diagram of a knot $K$, and let $F$ be a standard surface 
properly embedded in the knot exterior with at
most two meridional boundary components.
Then either $F$ is a trivial twice-punctured 2-sphere, 
or it admits a diagrammatic meridional compression.}

\noindent {\sl Proof.} 
We claim that we can find a simple closed $C$ of $F \cap S^2_+$
that intersects some bubble in at least two arcs, where two
of these arcs are part of the same saddle and have no
arc of $F \cap S^2_+$ between them. We can then 
find a simple closed curve in $F$
that runs from $C$ across the saddle back to $C$
and then over the disc that $C$ bounds above the
diagram. This curve bounds the required diagrammatic 
meridional compression disc.

If there is only one curve of $F \cap S^2_+$, then it
is either disjoint from the bubbles, in which case $F$ is a 
trivial twice-punctured 2-sphere, or it intersects some bubble
as claimed. Thus, we may assume that there are at least
two such curves. Consider an innermost one, $C$, bounding
a disc $I$ containing no other curves of $F \cap S^2_+$. 
By choosing this curve suitably, we can ensure that
$C$ has at most one point of intersection with $K$, since
$F$ has at most two boundary components.

Each time that $C$ runs over a bubble, the crossing either lies
in the inward direction or outward direction from $C$.
Note, however, that just because a crossing lies
in the inward direction from $C$, this does not guarantee
that the crossing itself lies in $I$, since
$C$ may return to the crossing several times. The 
hypothesis that the diagram is alternating ensures
that, as one runs along $C$, one meets crossings
alternately on the inward direction and outward direction of $C$.

Consider the arc components
of intersection between $I$ and the knot projection, ignoring
the components containing crossings. Suppose first that there
are at least two such arcs. Let $\alpha$ be an outermost
such arc in $I$, separating off a disc $E$. By choosing $\alpha$
appropriately, we can ensure that $E \cap C$
is disjoint from $K$. It therefore runs from a bubble
back to the same bubble. It meets this bubble once in the
inwards direction and once in the outwards direction.
Hence, by property (i) in the definition of a standard surface, 
it runs over at least
one other bubble in the inwards direction. 
Consider the curve of $F \cap S^2_+$ on the other side
of this crossing, connected via the saddle at the bubble.
This must again be part of $C$, by the assumption that
$C$ is innermost. Hence, we have the found the required
intersection with $C$ and a bubble.

Suppose now that there is precisely one arc $\alpha$ of
intersection between $I$ and the knot projection. If $C \cap K$
does not lie at an endpoint of $\alpha$, the above
argument works. Suppose therefore that $C \cap K$
does lie at an endpoint of $\alpha$. Note that $C$
must meet at least one other crossing, by (iii) in the definition
of a standard surface.
Hence, it meets a crossing in the inwards direction,
and the claim is then proved.

Suppose now that there is no arc of intersection
between $I$ and the knot projection. Now,
$C$ must meet at least three crossings, at least one of
which lies in the inward direction. Thus, as above,
the claim is proved. $\square$

\vskip 18pt
\centerline{\psfig{figure=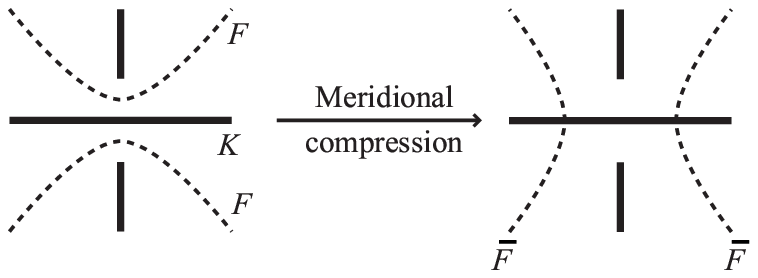}}
\vskip 18pt
\centerline{Figure 13.}

Let $\overline F$ be the result of $F$ 
immediately after the diagrammatic meridional compression. Then
$\overline F$ satisfies
all the conditions of a standard surface, except possibly (iii).
This may fail, since $\overline F$ may have one
or two `tubes' that run parallel to the knot for a time
before closing off with a disc that intersects the knot
in a single point. We now explain how to modify $\overline F$
by retracting these tubes, under the
assumption that the diagram is prime. Each
component of the resulting surface will be standard 
or a trivial 2-sphere.
For, if (iii) fails for a component of $\overline F$ that
is not a trivial 2-sphere, then there is then an
ambient isotopy which retracts this tube, reducing
the number of curves of $\overline F \cap S^2_+$ and leaving
the number of curves of $\overline F \cap S^2_-$ unchanged.
Clearly, (i) and (ii) still hold, and the surface
still intersects the bubbles in saddles, the faces in arcs
and the polyhedra in discs. Hence, eventually, each component is
a trivial 2-sphere or a standard surface.

\noindent {\bf Corollary 10.} {\sl Let $D$ be a prime
alternating diagram of a knot $K$. Then
\item{(i)} the exterior of $K$ contains no standard 2-spheres 
that intersect the knot in at most two points; 
\item{(ii)} any standard torus disjoint from the knot is
normally parallel to $\partial {\cal N}(K)$.

}

\noindent {\sl Proof.} Let $F$ be a standard 2-sphere disjoint from
the knot. Then by Lemma 9, $F$ admits a meridional compression.
The result is two 2-spheres, each of which intersects the knot
once, which is impossible. 

Now consider a standard 2-sphere $F$ intersecting $K$ in
two points. We will prove, by induction, on the number of its saddles
that such a 2-sphere cannot exist. It admits a meridional compression to
two such 2-spheres. Retract the tubes of these 2-spheres.
The result, inductively, cannot be standard 2-spheres.
Hence, they must be trivial. However, if we then
reconstruct $F$ by tubing these two trivial 2-spheres
together, the result is not standard.

Finally, let $F$ be a standard torus disjoint from the knot.
Find a diagrammatic meridional compression to a twice punctured
2-sphere. Retract its tubes. The result cannot be standard,
and hence must be trivial. Therefore, the original torus $F$
must have been normally parallel to $\partial {\cal N}(K)$.
$\square$

A further variation that we need to consider is 
alternating graphs as opposed to alternating knots.
In this case, we will need to analyse surfaces with
non-meridional boundary properly embedded in the graph
exterior. Again, there is a theory of normal
surfaces, following [6]. 
If one truncates the ideal vertices of the
ideal polyhedral decomposition, the resulting
truncated polyhedra have boundary that
can be identified with $S^2_\pm$. A properly embedded
surface is {\sl normal} if it
intersects each bubble in a collection of saddles,
intersects each face in a collection of arcs
and intersects the two truncated polyhedra in
a collection of discs. Also, the boundary $C$ of each such disc
must satisfy the following conditions:
\item{(i)} $C$ intersects each side
of each crossing in at most one arc;
\item{(ii)} $C$ does not lie entirely in $\partial {\cal N}(G)
\cap S^2_\pm$;
\item{(iii)} no arc of intersection between $C$ 
and a face $F$ has endpoints lying in the same crossing,
or in the same component of $\partial {\cal N}(G) \cap F$,
or in a crossing and a component of $\partial {\cal N}(G) \cap F$
that are adjacent;
\item{(iv)} $C$ intersects any component of 
${\cal N}(G) \cap S^2_\pm$ in at most one arc;
\item{(v)} any such arc cannot have endpoints in the same 
component of ${\cal N}(G) \cap F$ for any face $F$.

\vskip 18pt
\centerline{\caps 4. Alternating graphs with handlebody exteriors}
\vskip 6pt

The main goal of this section is to prove Theorem 3 below.

\noindent {\bf Theorem 3.} {\sl Let $D$ be a reduced
alternating diagram of a graph $G$ with a single
vertex. Then the exterior of $G$ is a handlebody
if and only if $D$ is one of the diagrams
shown in Figure 7.}

We will use the following lemma at a number of points.

\noindent {\bf Lemma 11.} {\sl Let $D$ be a reduced
alternating diagram of a graph $G$ with a single
vertex, such that the
exterior of $G$ is a handlebody $H$.  Let $C$ be a
simple closed curve in the diagram that intersects
the graph projection transversely in two points
disjoint from the crossings and the vertex. Then $C$ bounds
a disc in $D$ that is disjoint from the crossings
and the vertex.}

\noindent {\sl Proof.} Suppose, on the contrary,
that there is such a curve $C$, bounding a disc in
$D$ that is disjoint from the vertex but
contains at least one crossing. It bounds
two discs, one above the diagram, and one below.
The union of these is a 2-sphere, whose intersection
with $H$ is an annulus $A$. The core curve of $A$ is 
homologically non-trivial in $H$. Hence,
$A$ is incompressible in $H$. When a handlebody
is cut along a properly embedded orientable
incompressible surface, the result is a disjoint
union of handlebodies. Hence, the 1-string tangle
that $C$ bounds must be trivial. But this
contradicts Menasco's theorem [10], since it
is reduced, alternating and has at least one crossing. $\square$

\noindent {\sl Proof of Theorem 3.} We argued in the introduction
that the exteriors of the graphs in Figure 7
are handlebodies. Hence, we need only prove the converse.

Suppose now that the exterior of $G$ is a handlebody.
Its boundary therefore has a compression disc. Hence, by [6],
it has a compression disc $E$ in normal form. Note that
$E$ cannot have a meridional compression disc $E_1$.
For $\partial E_1$ would then bound a subdisc $E_2$ of $E$,
and $E_1 \cup E_2$ would be an embedded 2-sphere intersecting
$G$ in precisely one point, which is impossible.

Consider the intersection between $E$ and $S^2_- \cup S^2_+$.
This is a graph embedded in $E$. 
Its complimentary regions are saddles and normal discs,
where the former lie in the interior of $E$.
If this graph fails
to be connected, then pick a component that is outermost
in $E$. Let $E'$ be the subdisc of $E$ comprised of
this component and all faces adjacent to it. Suppose that the region
of $E'$ containing $\partial E' - \partial E$ 
lies below $S^2_-$, say. Then we will focus $E' \cap S^2_+$. This has the property that
any curve of $E' \cap S^2_+$ runs over $\partial {\cal N}(G)$ at most once.
Also, if $E' \cap S^2_+$ runs over a saddle, then
the curve of $E \cap S^2_+$ on the other side
of the crossing also lies in $E'$.

Let $N$ be the component of ${\cal N}(G) \cap S^2_+$
containing the vertex $v$.
If $N$ is not a disc, then one of the components $\beta$
of $G - v$ runs from $v$ back to $v$ without
passing underneath any crossings. Since $G$ is
alternating, it therefore runs through at most
one over-crossing. 
Suppose first that $\beta$ runs through no crossings.
Then the remainder of the diagram is a diagram of
the other component of $G - v$. By Lemma 11, this tangle has no crossings.
Then, $D$ is as in the leftmost diagram of Figure 7.
If $\beta$ runs through a single crossing,
then it divides the diagram into two 1-string tangles.
Each tangle has no crossings by Lemma 11. Hence, the diagram of $G$
has a single crossing, and therefore fails to be reduced.
We may therefore assume that $N$ is a disc.
Note that $\partial N$ runs through crossings eight times.

We perform a small ambient isotopy so that all
curves of $E' \cap S^2_+$ are disjoint from the
vertex $v$ of $G$. We can ensure that, for any curve
$C$ of $E' \cap S^2_+$, the disc of $S^2_+$ that $C$
bounds not containing $v$ contains at most four crossings
in $\partial N$.

Pick a curve $C$ of $E' \cap S^2_+$ innermost in
the diagram, where we define the innermost direction
so that the disc that $C$ bounds does not contain $v$.

Consider first the case where $C$ is disjoint from ${\cal N}(G)$.
Then it runs over an even number of crossings.
As in the proof of Lemma 7, the fact that the diagram
is alternating and that $C$ is innermost implies
that $E$ has a meridional compression disc, which
is impossible.

Thus $C$ must run over ${\cal N}(G)$. The arc
$C - {\cal N}(G)$ runs over at most one saddle. Otherwise,
we would either contradict the assumption that
$C$ is innermost or we would find a meridional
compression disc for $E$.

Suppose the component of $G \cap S^2_+$
that $C$ runs over does not contain the vertex of $G$.
Then it intersects $G$ in one of the ways shown in Figure 14.
In the first two cases, this part of $C$ intersects the projection of 
$G$ only once. Therefore, for parity reasons, $C - {\cal N}(G)$
must run over a crossing. In the third case of 
Figure 14, $C - {\cal N}(G)$ cannot run over
a crossing. Hence, $C$ bounds a diagram of a 1-string tangle
that must have no crossings. In the left diagram
of Figure 14, this contradicts the fact that $E$ is in normal form.
In the remaining two diagrams, we contradict the fact that $D$
is reduced and alternating.

\vskip 18pt
\centerline{\psfig{figure=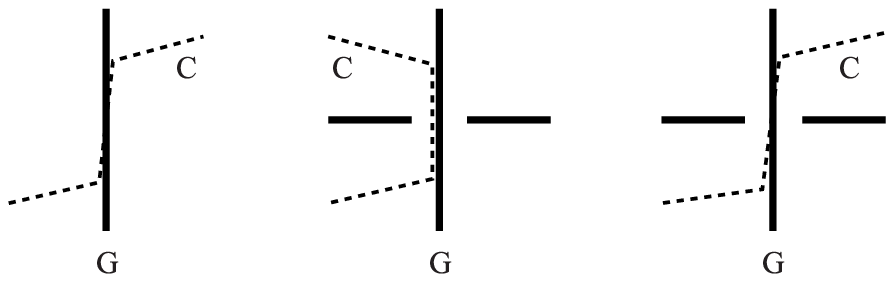}}
\vskip 18pt
\centerline{Figure 14.}

Thus, we are reduced to the case where $C$ runs
over the component $N$ of ${\cal N}(G) \cap S^2_+$ that
contains the vertex $v$. Let $\alpha$ be the arc ${\rm cl}(C - N)$.
Join the endpoints of $\alpha$ with a curve $\gamma$ in
$\partial N$  that is parallel in $N$ to $C \cap N$,
the parallelity region not containing $v$. Note
that $\gamma$ runs through at most
four crossings. The simple closed curve $\alpha \cup \gamma$
must run through an even number of crossings. If this number
is zero, then the disc is not in normal form. If this
number is two, then $\alpha \cup \gamma$ separates the diagram
into a connected sum. As before, $\alpha \cup  \gamma$ then bounds
a disc containing a single arc of the diagram and no crossings.
In each case, this implies that the diagram fails to be reduced
and alternating, or $E$ fails to be normal.
Hence, $\alpha \cup \gamma$ runs through precisely four
crossings. We list all the possibilities in Figure 15.
We have not included any possibilities that would
contradict the fact that the diagram is alternating. 

\vskip 18pt
\centerline{\psfig{figure=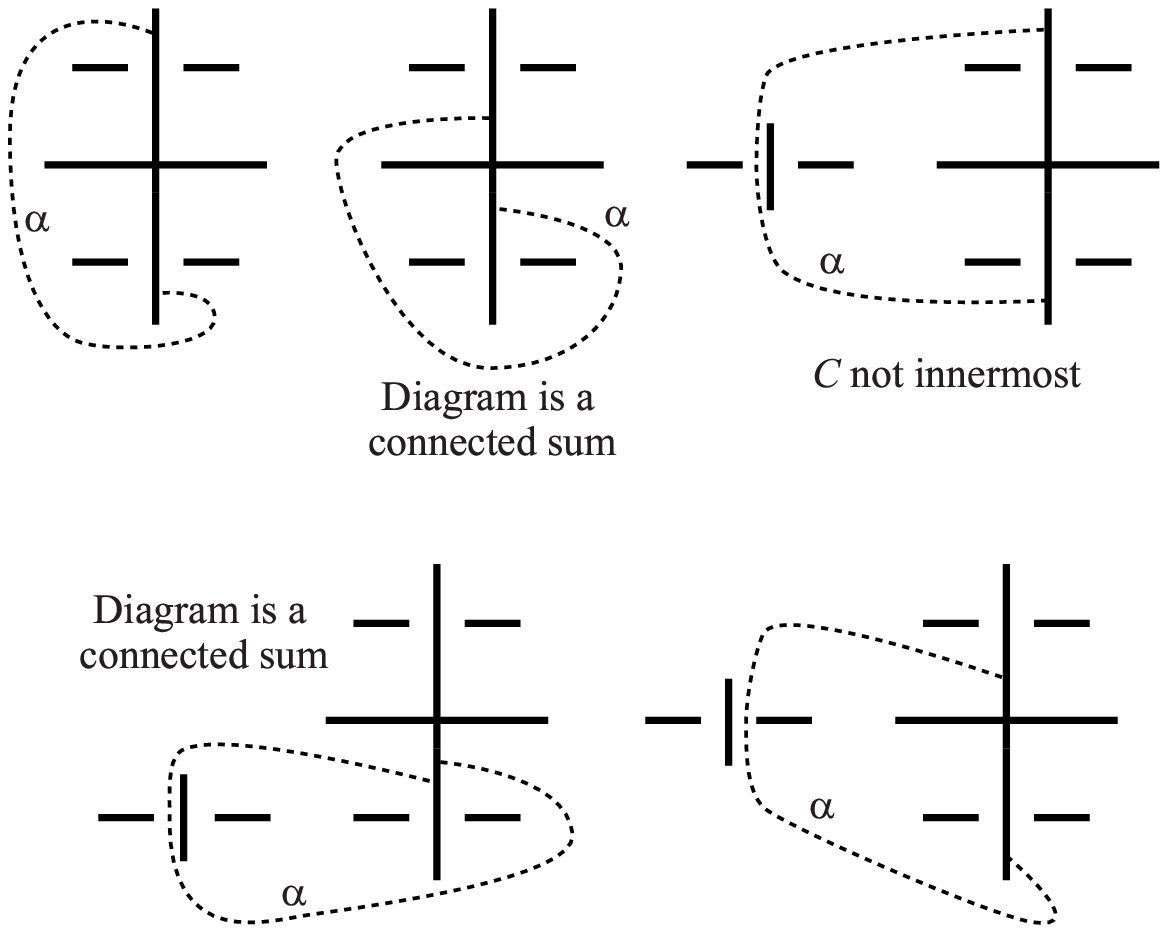,width=3in}}
\vskip 18pt
\centerline{Figure 15.}

It is clear that only two of these diagrams
may arise. In both cases, Figure 16 shows a way of removing
the vertex and replacing it with two arcs. 

\vskip 18pt
\centerline{\psfig{figure=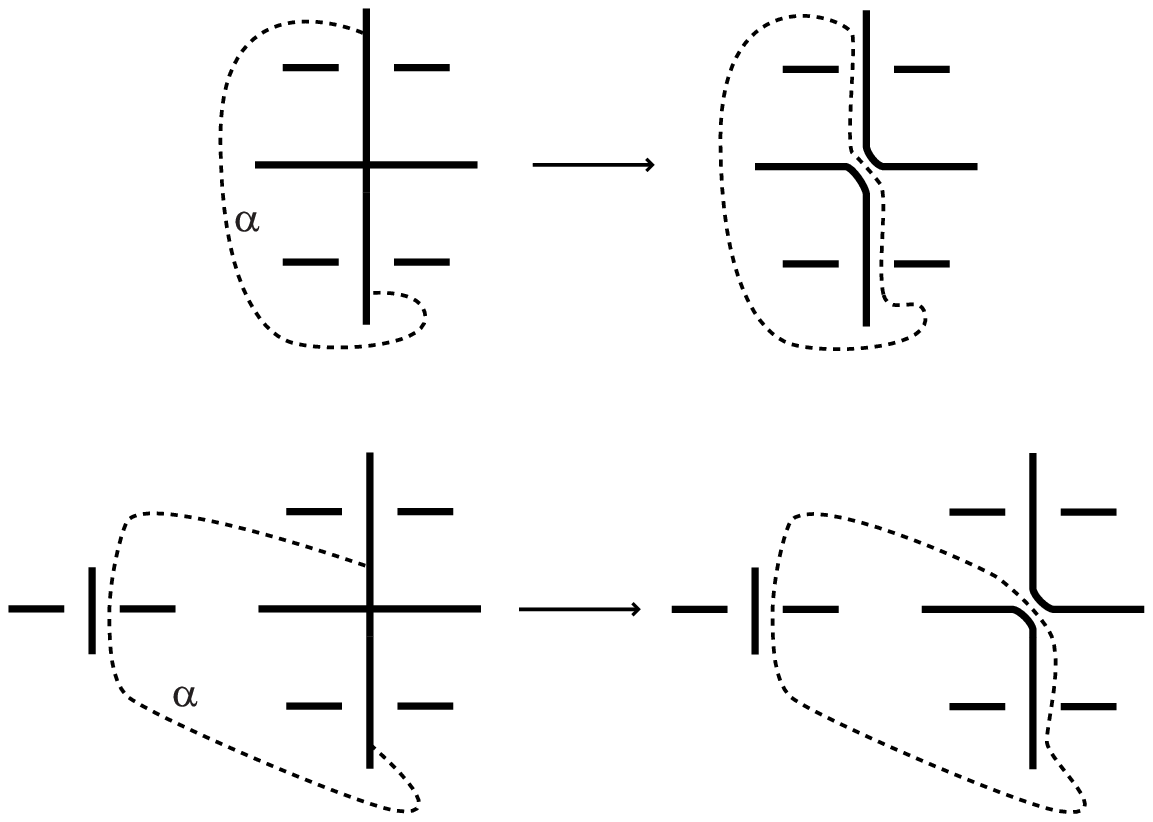,width=3in}}
\vskip 18pt
\centerline{Figure 16.}

The result is a composite alternating diagram of a knot or link $L$.
It has no trivial loops (namely, a loop that starts and ends
at the same crossing, with interior disjoint from the crossings).
Hence, $L$ is composite by [10]. This link $L$ has tunnel 
number one, since there is an obvious unknotting tunnel $t$
such that the exterior of $L \cup t$ is homeomorphic
to the exterior of $G$. A result of Gordon-Reid [4] asserts
that a tunnel number one knot or link can be composite only
if it has a Hopf link summand.
The Hopf link has a unique reduced alternating
diagram [11]. Therefore, in the diagram of $L$, one of
its components runs through exactly two crossings, forming
a meridian of the other component. Since $G$ is connected,
these two components are fused together at $v$. Thus,
a subset of $D$ is as shown in Figure 17. Note that this
is identical to the complement
of the grey discs in either the middle or right-hand
diagrams of Figure 7.

\vskip 18pt
\centerline{\psfig{figure=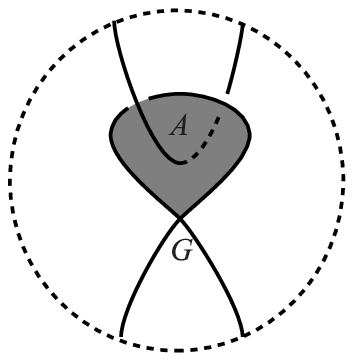}}
\vskip 18pt
\centerline{Figure 17.}

 We need to show that
the remainder of the diagram is a rational tangle.
The strings of this tangle must join up as shown in
Figure 7, since $G$ is connected.
Let $S$ be the 2-sphere bounding the subset of the
diagram in Figure 17.
This sphere divides $S^3$ into two 3-balls, $B_1$ and
$B_2$, where we say that $B_2$ is the one containing the
vertex $v$.

We need to show that the tangle $G \cap B_1$ is rational. 
Let $A$ be the annulus shown in Figure 17 properly
embedded in the exterior of $G$. Note that
there is a homeomorphism $h$ from $S^3 - {\rm int}({\cal N}(G \cup A))$
to $B_1 - {\rm int}({\cal N}(G \cap B_1))$, taking
the two copies of $A$ in $\partial {\cal N}(G \cup A)$ to
the two annuli ${\rm cl}(\partial {\cal N}(G \cap B_1) - \partial B_1)$.
This annulus $A$ is incompressible in the complement of
$G$. But the exterior $H$ of $G$ is a handlebody. Hence,
$A$ admits a $\partial$-compression in $H$. The image under
$h$ of this $\partial$-compression disc is a disc $P$
embedded in $B_1 - G$. Note that $P \cap \partial B_1$ is a
single arc, and the remainder of $\partial P$ runs once
along one component of $G \cap B_1$. Hence, $P$ specifies
a parallelity disc between this component of $G \cap B_1$
and an arc in $\partial B_1$. The remaining arc
of $G \cap B_1$ lies in the complement of $P$, which is
a 3-ball. It must be a trivial 1-string tangle in this
3-ball. Hence, the intersection of $B_1$ with $G$ is a
rational tangle. This proves that $G$ is as shown in Figure 7.
$\square$

\vskip 18pt
\centerline{\caps 5. Almost normal surfaces in alternating knot
complements}
\vskip 6pt

Our goal in this section is to prove Theorem 2, and hence
complete the proof of Theorem 1.

\noindent {\bf Theorem 2.} {\sl Let $D$ be a reduced alternating diagram
for a knot $K$. Then any unknotting tunnel for $K$ is
isotopic to an unknotting tunnel that is an embedded arc 
in some region of the diagram.}

Note first that the diagram $D$ is prime. Otherwise, $K$ is a non-trivial
connected sum [10], which is impossible by Gordon-Reid's theorem [4].

\noindent {\bf Step 1.} The Heegaard surface is strongly irreducible.

We will denote the genus two Heegaard surface by $F$.
Let $M$ be $S^3 - {\rm int}({\cal N}(K))$. Let $H$ and $P$ be the closures
of the two components of $M - F$, where $H$ is the handlebody
and $P$ is the compression body containing $\partial {\cal N}(K)$
in its boundary. Note first that we may assume that $F$ is
irreducible. For otherwise, the Heegaard splitting is the connected
sum of a genus one splitting of $S^3$ and a genus one splitting
for the knot exterior. By Waldhausen's theorem [20], the splitting
is stabilised. Thus, the unknotting tunnel is isotopic to a planar
arc, which proves the theorem. 

Suppose that $F$ is weakly reducible, via disjoint compressing
discs $D_1$ and $D_2$ in $H$ and $P$ respectively. If
we cut $P$ along $D_2$, the result is 
a copy of $T^2 \times I$ and possibly also a solid torus $V$.
The curve $\partial D_1$ is disjoint from $D_2$.
Since the splitting is irreducible, it does not bound a disc in $P$. 
If it lies in $T^2 \times I$, it extends to a
compression disc for $\partial {\cal N}(K)$. Then $K$ is
the unknot, and by [2], the splitting is reducible, contrary to
assumption. If $\partial D_1$ lies in the solid torus,
then we obtain a lens space
summand for the knot exterior. This lens space
must be $S^3$ and hence $\partial D_1$ has winding number one round
$V$. This implies that the Heegaard splitting
is reducible, which is a contradiction.
This proves that $F$ is strongly irreducible.

\noindent {\bf Step 2.} Placing $F$ into almost normal form.

Pick some subdivision of the faces of the ideal polyhedral decomposition
of the knot complement, without introducing any vertices,
so that each face is bigon or triangle.
By Corollary 10, this contains no properly embedded
standard 2-spheres. In particular, it contains no 2-spheres
that are normal to one side. So, by Theorem 5, we may ambient isotope $F$ into
almost normal form in this ideal polyhedral decomposition.

Consider first the case where $F$ has an almost normal tubed
piece. Compress this tube. The result is a surface that
is normal to one side. By Corollary 10, it contains no 2-sphere
components. Hence, it is one or two tori. In the latter case, one
of these tori would be compressible. Each torus is
standard, and hence, by Corollary 10, is normally parallel
to $\partial {\cal N}(K)$.
In particular, it is incompressible. Thus, the compressed surface
is a single torus $\overline F$ that is normally parallel to $\partial {\cal N}(K)$. 
Up to isotopy, the unknotting tunnel runs
from $\partial {\cal N}(K)$ through the tube and then to
$\partial {\cal N}(K)$, respecting the product structure
in the tube and the product structure on the parallelity
region between $\partial {\cal N}(K)$ and $\overline F$. This arc
can therefore be ambient isotoped into a face of the
ideal polyhedral decomposition. This is then an embedded
arc in some region of the diagram. This proves Theorem 2 in
this case.

We therefore assume now that the almost normal surface $F$ 
contains no almost normal tubed piece.

\noindent {\bf Step 3.} Meridionally compressing $F$.

By Lemma 9, $F$ admits a diagrammatic meridional compression
to a twice-punctured torus $T$. Retract the tubes of $T$
to place it in standard form. Now apply Lemma 9 again
to find another meridional compression disc.
There are two possibilities for the
resulting surface: {\sl either} a 2-sphere 
intersecting $K$ in four points {\sl or} a
torus intersecting $K$ in two points and
a 2-sphere intersecting $K$ twice.
We claim that the latter case cannot arise.
For the twice-punctured 2-sphere retracts to
a trivial 2-sphere by Corollary 10. But, then
reconstructing $T$, we see that it could not
have been standard. Hence, the meridional compression must
yield a 2-sphere $S$ intersecting $K$ in four points.

\noindent {\bf Step 4.} How the tubes can be nested.

When we performed the first meridional compression to $F$, this
created two points of intersection between $T$ and $K$. Let
$\alpha_1$ be the sub-arc of $K$ that runs between these two points,
and which lies in a regular neighbourhood of the compression
disc. When we retracted the tubes of $T$, this expands
$\alpha_1$, but it remains an embedded sub-arc of $K$.
When we perform the second meridional compression, we get
a similar arc $\alpha_2$ running along $K$ between two points of
$S \cap K$. The surface $F$ is obtained from $S$ by
removing the four discs $S \cap {\cal N}(K)$
and attaching a tube (that is, an annulus) that runs along
$\alpha_2$ and then another that runs along $\alpha_1$.
There are two cases to consider: either
$\alpha_1$ and $\alpha_2$ are disjoint, or $\alpha_2$
is a subset of $\alpha_1$. In the latter case, we say
that the tubes are nested, whereas in the former case,
they are not.

\noindent {\bf Step 5.} $S$ bounds a rational tangle
on at least one side.

Consider first the case where the tubes are not nested.
Since $F$ is a Heegaard surface, $K$ is parallel to a
simple closed curve on $F$, via an annulus.
We may assume that the meridional compression
discs each intersect this annulus in a single arc, and
hence, when we cut along these discs, we obtain two
discs $E_1$ and $E_2$, where $E_i \cap K$ is an
arc in $\partial E_i$, and $E_i \cap S$ is the
remainder of $\partial E_i$. Hence, in this case,
$S$ bounds a rational tangle on one side. We call
the 3-ball on this side $B_1$. 

Consider now the case where the tubes are nested.
Let $B_1$ and $B_2$ be the 3-balls on each side of $S$.
The Heegaard surface $F$ is obtained from $S$ in two stages:
\item{(i)} attaching a tube $T_2$ running along $K$, 
joining two punctures of $S$; suppose that this tube lies 
in the 3-ball $B_2$ bounded by $S$, and then
\item{(ii)} by attaching another tube $T_1$; this time $T_1$ runs between
the two other punctures of $S$; it runs into $B_1$, back to $S$,
through the tube formed in stage (i), and then back through $B_1$.

The component of $S^3 - {\rm
int}({\cal N}(F))$ not containing $K$ is a handlebody $H$.
This component consists of $B_1 - {\rm int}({\cal N}(K \cap B_1))$
together with the space between the two tubes, which is
a copy of $A \times I$, where $A$ is an annulus and
where $(A \times I) \cap B_1 = A \times \partial I$.
Then $A \times \{ {1 \over 2} \}$ 
is an incompressible annulus properly embedded in a handlebody.
Such an annulus must have a boundary-compression disc $E$.
We may assume that $E$ intersects $A \times I$ as the
product of a properly embedded essential arc in $A$
and either $[0, {1 \over 2}]$ or $[{1 \over 2}, 1]$.
We may also assume that $E$ intersects the remainder of $T_1$
in a single arc. By extending $E - (A \times (0,1))$ a
little, we obtain a disc $E'$ 
embedded in $B_1$, such that $\partial E'$ is the union
of an arc of $K \cap B_1$ and an arc in $\partial B_1$.
The other arc of $K \cap B_1$ lies in $B_1 - {\rm int}({\cal N}(E'))$,
which is a 3-ball. This arc must also be trivial,
since $H - {\rm int}({\cal N}(A))$ is a handlebody.
Hence, $(B_1, B_1 \cap K)$ is a rational tangle.

\noindent {\bf Step 6.} The possible positions for $S$.

We retract the tubes of $S$ to a standard surface. 
In [10], Menasco analysed
in some detail the possible arrangements for a
four-times punctured 2-sphere in normal form.
Although $S$ is not necessarily normal, Menasco's
arguments still hold here. We claim that $S$
intersects $S^2_+$ and $S^2_-$ as in one of the possibilities
of Figure 18.

We first claim that $S$ has no diagrammatic meridional compression disc. 
For, the boundary of this disc would have linking number one with $K$. However,
the only curves on $S - K$ with this property are
parallel to one of the curves of $S \cap \partial {\cal N}(K)$.
Hence, a diagrammatic meridional compression would yield another
four-times punctured 2-sphere and a twice-punctured
2-sphere. We argued in Step 3 that there can be no
such twice-punctured 2-sphere.

\vskip 18pt
\centerline{\psfig{figure=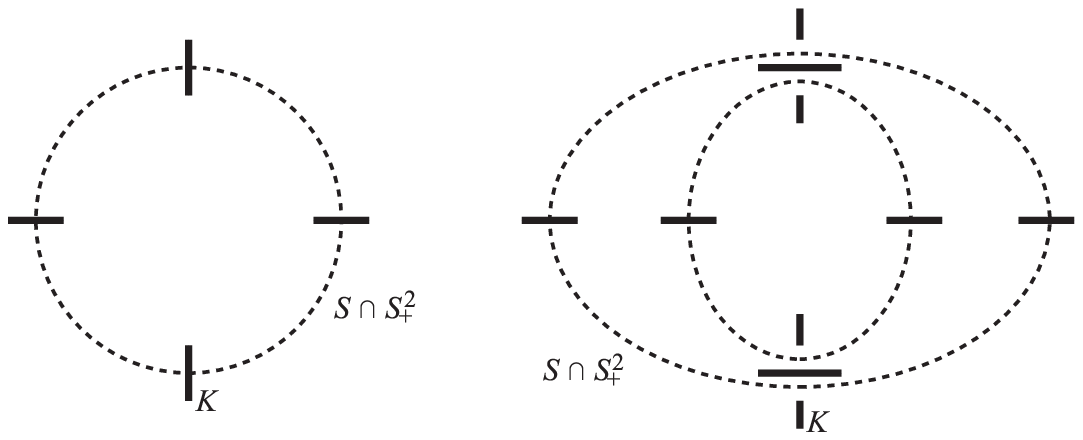}}
\vskip 18pt
\centerline{Figure 18.}

Applying the arguments of Lemma 9 and the fact that
$S$ has no diagrammatic meridional compression, we deduce that
each innermost curve
of $S \cap S^2_\pm$ must intersect $K$ at least twice.
If there is a single curve of $S^2_+ \cap S$, then it 
must meet $K$ four times and can meet no saddles.
It is then as in shown in the left-hand diagram of Figure 18,
which we term the simple arrangement for $S$. If there is more than
one curve of $S^2_+ \cap S$, then there at least two that
are innermost. Each of these curves has two points of intersection
with $K$. As in the proof of Lemma 9, they cannot meet a saddle in the inwards
direction. Hence, each meets two saddles, the saddles alternating with
points of $K$. Since this holds in both $S^2_+$ and $S^2_-$, it
is not hard to see that the only possibility is as in the
right-hand diagram of Figure 18, which we term the
complex arrangement for $S$.

Our aim is to reduce the complex arrangement to the
simple arrangement. We will show that, in the complex
case, there is an ambient isotopy, leaving $K$ invariant
and introducing no new point of $S \cap K$, taking
$S$ to a simple arrangement.
So, suppose that $S$ is as in the right-hand diagram
of Figure 18. There is an obvious vertical disc $E$ with 
$E \cap S = \partial E$,
which intersects $K$ in two points. The boundary of this disc
runs from the top crossing, along the inner disc of $S - S^2_+$
above $S^2_+$ to the next crossing, under the saddle, then
over the top disc of $S - S^2_+$ back to the top crossing.
There is a similar disc under $S^2_-$. We choose $E$ 
so that it lies in $B_1$, the 3-ball containing the rational
tangle. Note that $\partial E$ separates the four
points of $S \cap K$ into two pairs.

Let $P$ be two disjoint discs embedded in $B_1$, so that $\partial P$
contains the two arcs $B_1 \cap K$, one in each component
of $P$, and so that the remainder of $\partial P$ is
$P \cap \partial B_1$. Such discs $P$ exist because
$(B_1, B_1 \cap K)$ is a rational tangle. We may isotope
$P$ so that it intersects $E$ in simple closed curves
and embedded arcs. The six possible curve or arc
types of $E \cap P$ in $E$ are shown in Figure 19.
Note that a curve of $E \cap P$ encircling a single
point of $E \cap K$ is ruled out, since the arc of
$E \cap P$ emanating from this point must have another
endpoint.

\vfill\eject
\centerline{\psfig{figure=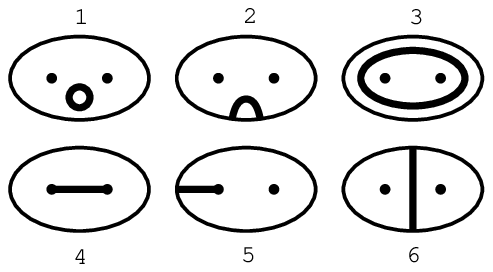}}
\vskip 18pt
\centerline{Figure 19.}

If there is a type 1 curve, pass to an innermost one
in $E$. This bounds a disc in $P$ and a disc in $E$.
Alter $P$, replacing the former disc with the latter,
and then perform a small ambient isotopy to remove this
curve of $E \cap P$. In this way, we remove all type 1
curves. Similarly, by dealing with outermost type 2 curves
in $E$, we may remove all type 2 curves. There are now
two possibilities: $E$ contains a type 4 curve and perhaps
some type 3 curves; {\sl or} $E$ contains two type 5 curves
and perhaps some type 6 curves.

We start with the first case. Suppose that there is a
type 3 curve. An innermost type 3 curve in $P$ bounds a
disc whose interior is disjoint from $E$. The union of
this disc and $E$ divides $B_1$ into three balls.
The 3-ball that intersects $E$ in an annulus
contains a single arc of $B_1 \cap K$, which
must be a trivial tangle. Similarly, the 3-ball disjoint from 
$\partial B_1$ intersects $K$ in a single arc $\alpha$, which again must 
be a trivial tangle. Hence, there is a disc embedded in this
ball, whose boundary is $\alpha$ and a single
type 4 curve. This type 4 curve is isotopic to $\alpha$,
leaving its boundary fixed. Hence,
the final 3-ball $B'_1$ (the one containing all of $E$ in
its boundary) intersects $K$ in a 2-string tangle,
so that the union of these two strings, the type
4 curve and an arc in $\partial B_1 - E$ bounds a
disc in $B'_1$. (See Figure 20.) Hence, we may reconstruct $P$ so that it
intersects $E$ in a single type 4 curve. We may
therefore assume that $P$ contains no type 3 curves.

\vskip 18pt
\centerline{\psfig{figure=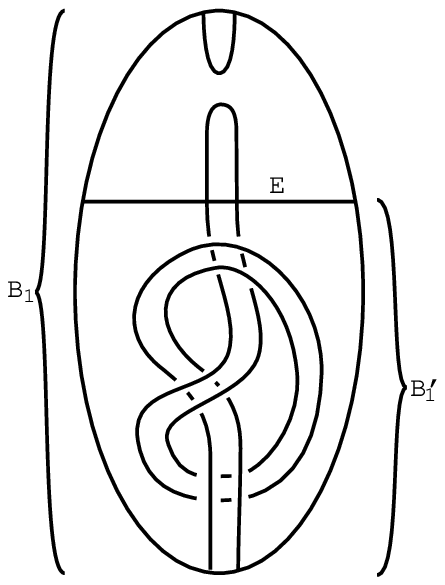}}
\vskip 18pt
\centerline{Figure 20.}

We claim that the `satellite' tangle in $B'_1$ is
trivial. In other words, pair $(B'_1, B'_1 \cap K)$ is
homeomorphic to the product of an interval and the
pair $(E, E \cap K)$.
The tangle in $B'_1$ has a diagram which is a subset
of the alternating diagram $D$. However, a non-trivial
satellite tangle cannot have an alternating diagram.
Otherwise, we could extend the tangle to an alternating
diagram of a prime non-trivial satellite knot, contradicting
Menasco's theorem [10]. This proves the claim. We may now use
the product structure on $B'_1$ to isotope $S$, as
required, across $B'_1$, so that afterwards $S$
has a simple arrangement in the diagram.

We now deal with the case where $P$ intersects $E$
in two type 5 curves and possibly some type 6 curves. If there
is a type 6 curve, consider one outermost in $P$.
This separates off a subdisc $P'$ of $P$ that is
disjoint from $K$. It lies 
in one component $B'_1$ of ${\rm cl}(B_1 - E)$.
It separates the two arcs of $B'_1 \cap K$.
The intersection between $K$ and each component of
$B'_1 - P'$ is a trivial 1-string tangle.
Since $\partial P'$ runs over $E$ in a single arc,
we deduce again that the pair $(B'_1, B'_1 \cap K)$ is
homeomorphic to the product of an interval and the
pair $(E, E \cap K)$. Similarly, if there are no
type 6 curves, the closures of both components of $B - E$
have such a product structure. Hence, again,
there is an ambient isotopy,
leaving $K$ invariant and introducing no new points of
$S \cap K$, taking $S$ to the 2-sphere $(S - \partial B'_1) \cup E$,
which has a simple arrangement with the knot diagram.

Hence, we may assume that $S \cap S^2_\pm$ is
a single simple closed curve containing all four points of
$S \cap K$.

In the case where the tubes of $F$ are not nested, the proof
of Theorem 2 is now complete. Up to isotopy, the
unknotting tunnel lies in the rational tangle $B_1$
as shown in Figure 21. This is a planar arc in the diagram.

\vskip 18pt
\centerline{\psfig{figure=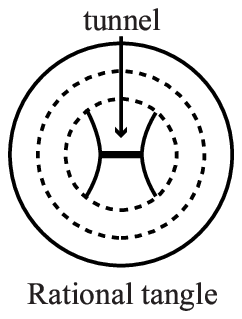}}
\vskip 18pt
\centerline{Figure 21.}

\noindent {\bf Step 7.} The case where the tubes are nested.

Recall from Step 5 that $S$ divides $S^3$ into $B_1$ and $B_2$.
The tangle $(B_1, B_1 \cap K)$ is rational. 
We now need to analyse the other tangle $(B_2, B_2 \cap K)$.
We know that $F$ is a Heegaard surface for $K$. 
Hence, the component $H'$ of $S^3 - {\rm int}({\cal N}(F))$
containing $K$ is obtained by attaching a 1-handle
to a solid torus in which $K$ is the core curve.
There is therefore an annulus $A'$ embedded in this
handlebody $H'$, with $K$ as one boundary component
and the other boundary component lying in $\partial H'$.
If $D_1$ is the first meridional compression disc for
$F$, we may assume that $D_1 \cap A'$ is a single arc
running from $\partial D_1$ to $D_1 \cap K$.
Hence, $D' = A' - {\rm int}({\cal N}(D_1))$
is a disc lying on one side of $T$, such that
$D' \cap (T \cup K) = \partial D'$, with
$\partial D'$ comprising an arc in $K$ and an arc in $T$.
Note also that $T$ has a compression
disc disjoint from $K$ on the $D'$ side of $T$. Therefore, the component
of $S^3 - {\rm int}({\cal N}(T))$ containing $D'$
is a solid torus $V$, and $K \cap V$ is a curve $\alpha_3$
parallel to an arc in $\partial V$, via the disc $D'$. This solid torus
is $B_2 - {\rm int} ({\cal N}(\alpha_2))$. We will exhibit
a planar arc $t$ with endpoints in $\alpha_3$ such that
$\alpha_3 \cup t$ is the union of a core curve in $V$
and two vertical arcs. This $t$ will therefore be
isotopic to the original unknotting tunnel.

Now, $B_2 - {\rm int}({\cal N}(\alpha_2 \cup \alpha_3))$ is
a handlebody. We may view the subset of the diagram
that specifies $B_2$ as an alternating
diagram for this handlebody, where the outside
of $B_2$ is a single vertex. This diagram
need not be reduced, since it may be as in the
right-hand diagram of Figure 5. But nevertheless,
we deduce from Theorem 3 that the diagram for $B_2$ is
either rational or as shown in Figure 22, where the arcs
in $R$ join the boundary components so that this tangle
has two strings. As usual, the restriction of the diagram
to each of the concentric annuli has a single crossing and
contains four arcs, each running between the two boundary components of the
annulus.

\vskip 18pt
\centerline{\psfig{figure=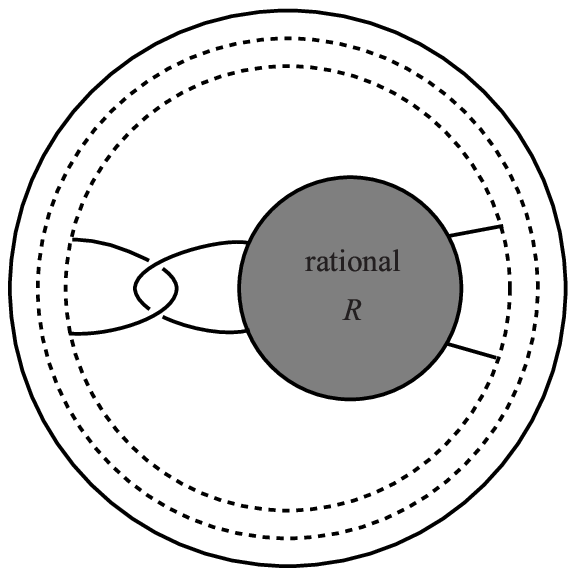,width=1.7in}}
\vskip 18pt
\centerline{Figure 22.}

Consider first the case where $B_2$ is rational. 
The diagrams for $B_1$ and $B_2$ are as shown in Figure 2.
Take $t$ in $B_2$ to be the tunnel shown
in the left diagram of Figure 23, where the endpoints
of $t$ lie in $\alpha_3$. Then $\alpha_3 \cup t$ is the 
union of a core curve in $V$ and two vertical arcs, as
required. However, $t$ is not yet planar.
Consider the two annuli closest to the
disc containing this tunnel, as shown in Figure 23. In each case, there is
an ambient isotopy that makes this tunnel planar. This proves
the theorem in this case.

\vskip 18pt
\centerline{\psfig{figure=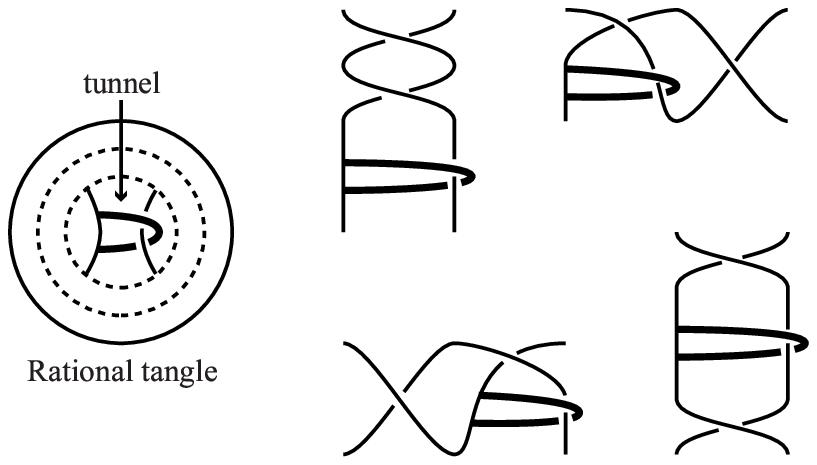}}
\vskip 18pt
\centerline{Figure 23.}

Now suppose that $B_2$ is as shown in Figure 22.
Note that we may assume that $R$ contains
at least two crossings, for otherwise the tangle in Figure 22 is
rational.
There are two cases: when the left-hand arc in Figure 22 
(that runs through only two crossings not in the outer annuli) 
is $\alpha_3$, and when it is $\alpha_2$. When it is
$\alpha_3$, we know that 
the tangle $R$ is integral, as shown in Figure 24,
since $\alpha_2$ must be a trivial 1-string tangle in $B_2$,
as $V = B_2 - {\rm int}({\cal N}(\alpha_2))$ is a solid torus.
The arc $t$ shown in Figure 24 is the required 
unknotting tunnel.

\vskip 18pt
\centerline{\psfig{figure=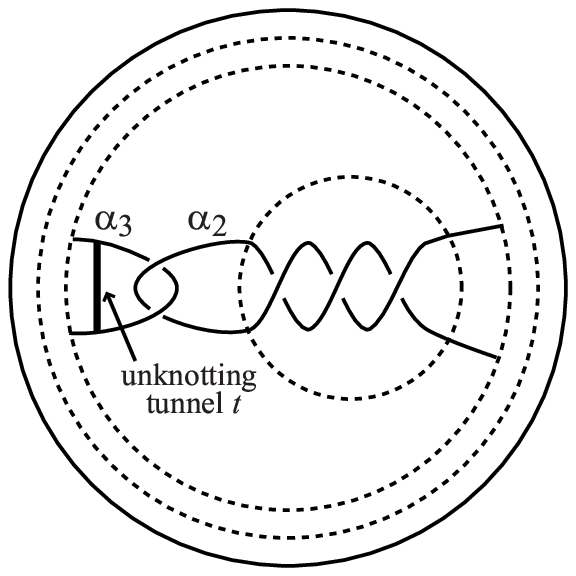,width=1.5in}}
\vskip 18pt
\centerline{Figure 24.}

When the left-hand arc of Figure 22 is $\alpha_2$, we 
take $t$ to be the planar tunnel shown in Figure 25. We must
show that $t \cup \alpha_3$ is the union of
a core curve of $V$ and two vertical arcs.
We will ambient isotope $\alpha_3 \cup t$ in the solid torus $V$
until this claim is evident. In fact, $t$ will remain
fixed.

\vskip 18pt
\centerline{\psfig{figure=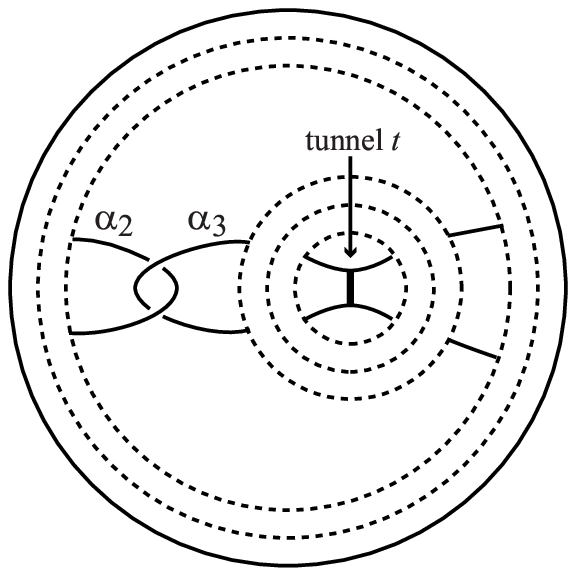,width=1.5in}}
\vskip 18pt
\centerline{Figure 25.}

By applying a homeomorphism to $B_2 - {\rm int}({\cal N}(\alpha_2))$,
we may remove the crossings in the outer annuli.
The tangle $R$ is decomposed into annular diagrams, each
containing a single crossing, surrounding a 2-string diagram
with no crossings. Consider the crossing of $R$ lying in the
outermost annulus $A_1$.
If it joins the right two points of the annulus,
then it may be removed, since it can
be viewed as a crossing in an outer annulus of the diagram for $B_2$. 
If it joins the left two points, then we may flype $R$ so that it
lies between the right two points, and then remove it. The cases where
the crossing joins the top two points, and where it joins the bottom
two points are symmetric. Hence we suppose the former.
Now consider the second outermost annulus $A_2$. If there is
none, then $R$ has only one crossing, and we are done.
Flype the subset of $R$ inside $A_2$, if necessary,
so that the crossing in $A_2$ joins the top two points,
or the left two points. There are thus two cases, which
are shown in Figure 26. There, an isotopy of $\alpha_3$
in $V$ is shown which removes these two crossings,
but leaves the remaining crossings unchanged.

\vskip 18pt
\centerline{\psfig{figure=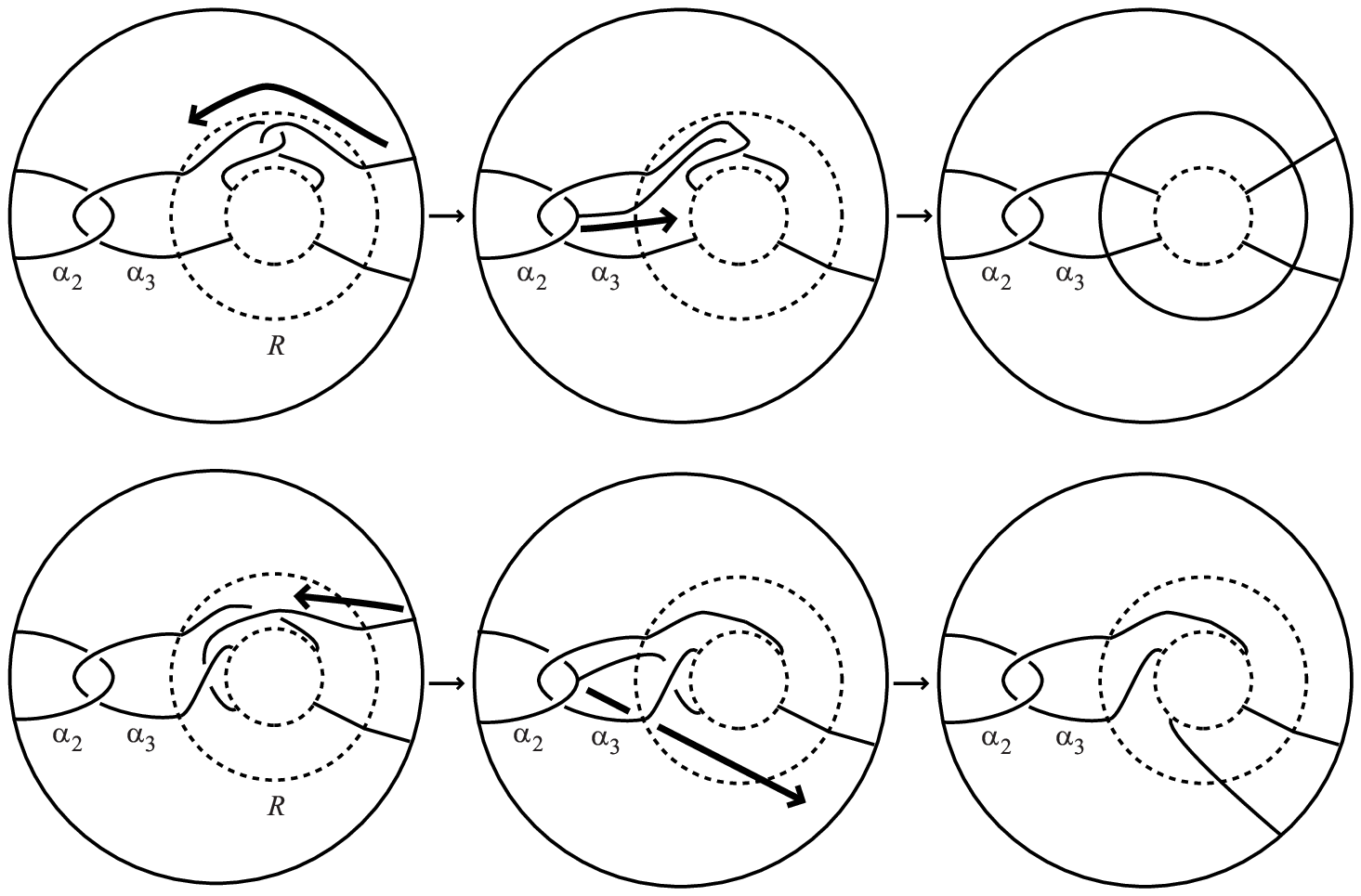,width=4in}}
\vskip 18pt
\centerline{Figure 26.}

It is clear therefore that we can isotope $\alpha_3$
until $R$ has at most one crossing. But, in this case,
$\alpha_3 \cup t$ is clearly the union of a core curve
and two vertical arcs.

This completes the proof of Theorem 2.

\vfill\eject
\centerline{\caps References}
\vskip 6pt

\item{1.} G. BURDE and H. ZIESCHANG, {\sl Knots}, de Gruyter (1985)
\item{2.} A. CASSON and C. GORDON, {\sl Reducing Heegaard splittings},
Topology and its Appl. {\bf 27} (1987) 275--283.
\item{3.} R. CROWELL, {\sl Genus of alternating link types},
Ann. of Math. (2) {\bf 69} (1959) 258--275.
\item{4.} C. GORDON and A. REID, {\sl Tangle decompositions of tunnel number 
one knots and links,} J. Knot Theory Ramifications {\bf 4} (1995) 389--409.
\item{5.} J. HEMPEL, {\sl 3-Manifolds}, Ann. of Math. Studies, No. 86,
Princeton Univ. Press, Princeton, N. J. (1976)
\item{6.} W. JACO and U. OERTEL, {\sl An algorithm to decide if a 3-manifold
is Haken}, Topology {\bf 23} (1984) 195--209.
\item{7.} L. KAUFFMAN, {\sl State models and the Jones polynomial,}
Topology {\bf 26} (1987) 395--407. 
\item{8.} T. KOBAYASHI, {\sl Classification of unknotting tunnels for two 
bridge knots.} Proceedings of the Kirbyfest (Berkeley, CA, 1998), 259--290 
Geom. Topol. Monogr., 2, Geom. Topol., Coventry (1999).
\item{9.} W. MENASCO, {\sl Polyhedra representation of link complements},
Low- \hfill\break
dimensional Topology, Contemp. Math {\bf 20}, Amer. Math. Soc.
(1983) 305--325.
\item{10.} W. MENASCO, {\sl Closed incompressible surfaces in
alternating knot and link complements}, Topology {\bf 23}
(1984) 37--44.
\item{11.} W. MENASCO and M. THISTLETHWAITE, {\sl Surfaces with
boundary in alternating knot exteriors}, J. Reine Angew. Math. 
{\bf 426} (1992) 47--65.
\item{12.} K. MURASUGI, {\sl On the genus of the alternating knot. I, II,}
J. Math. Soc. Japan {\bf 10} (1958) 94--105, 235--248.
\item{13.} K. MURASUGI, {\sl Jones polynomials and classical conjectures in 
knot theory,} Topology {\bf 26} (1987) 187--194. 
\item{14.} J. H. RUBINSTEIN, {\sl Polyhedral minimal surfaces,
Heegaard splittings and decision problems for 3-dimensional
manifolds}, Proceedings of the Georgia Topology Conference,
AMS/IP Stud. Adv. Math, vol. 21, Amer. Math Soc. (1997) 1--20.
\item{15.} K. SHIMOKAWA, {\sl On tunnel number one alternating
knots and links}, J. Math. Sci. Univ. Tokyo {\bf 5} (1998) 547--560.
\item{16.} M. STOCKING, {\sl Almost normal surfaces in 3-manifolds},
Trans. Amer. Math. Soc. {\bf 352} (2000) 171--207.
\item{17.} M. THISTLETHWAITE, {\sl
A spanning tree expansion of the Jones polynomial,}
Topology {\bf 26} (1987) 297--309. 
\item{18.} M. THISTLETHWAITE, {\sl On the algebraic part of an
alternating link}, Pacific J. Math. {\bf 151} (1991) 317--333.
\item{19.} A. THOMPSON, {\sl Thin position and the recognition
problem for the 3-sphere}, Math. Res. Lett. {\bf 1} (1994) 613--630.
\item{20.} F. WALDHAUSEN, {\sl Heegaard-Zerlegungen der $3$-Sph\"are.} 
Topology {\bf 7} (1968) 195--203.

\vskip 12pt
\+ Mathematical Institute,\cr
\+ Oxford University, \cr
\+ 24-29 St Giles',\cr
\+ Oxford OX1 3LB, \cr
\+ England.\cr

\end